\documentclass[reqno]{gokova}
\usepackage{epsfig,graphicx}

\begin{document}
\def\E{\ifmmode{\mathbb E}\else{$\mathbb E$}\fi} %natural numbers
\def\N{\ifmmode{\mathbb N}\else{$\mathbb N$}\fi} %natural numbers
\def\R{\ifmmode{\mathbb R}\else{$\mathbb R$}\fi} %real numbers
\def\Q{\ifmmode{\mathbb Q}\else{$\mathbb Q$}\fi} %rational numbers
\def\C{\ifmmode{\mathbb C}\else{$\mathbb C$}\fi} %complex numbers
\def\H{\ifmmode{\mathbb H}\else{$\mathbb H$}\fi} %complex numbers
\def\Z{\ifmmode{\mathbb Z}\else{$\mathbb Z$}\fi} %integers
\def\P{\ifmmode{\mathbb P}\else{$\mathbb P$}\fi} %real numbers
\def\T{\ifmmode{\mathbb T}\else{$\mathbb T$}\fi} %real numbers
\def\SS{\ifmmode{\mathbb S}\else{$\mathbb S$}\fi} %real numbers
\def\DD{\ifmmode{\mathbb D}\else{$\mathbb D$}\fi} %real numbers

\renewcommand{\a}{\alpha}
\renewcommand{\b}{\beta}
\renewcommand{\d}{\delta}
\newcommand{\D}{\Delta}
\newcommand{\e}{\varepsilon}
\newcommand{\g}{\gamma}
\newcommand{\G}{\Gamma}
\newcommand{\la}{\lambda}
\newcommand{\La}{\Lambda}
\newcommand{\n}{\nabla}
\newcommand{\var}{\varphi}
\newcommand{\s}{\sigma}
\newcommand{\Sig}{\Sigma}
\renewcommand{\t}{\tau}
\renewcommand{\th}{\theta}
\renewcommand{\O}{\Omega}
\renewcommand{\o}{\omega}
\newcommand{\z}{\zeta}

\newcommand{\ben}{\begin{enumerate}}
\newcommand{\een}{\end{enumerate}}
\newcommand{\be}{\begin{equation}}
\newcommand{\ee}{\end{equation}}
\newcommand{\bea}{\begin{eqnarray}}
\newcommand{\eea}{\end{eqnarray}}
\newcommand{\bc}{\begin{center}}
\newcommand{\ec}{\end{center}}

\newcommand{\IR}{\mbox{I \hspace{-0.2cm}R}}
\newcommand{\IN}{\mbox{I \hspace{-0.2cm}N}}

\newtheorem{thm}{Theorem}[section]
\newtheorem{cor}[thm]{Corollary}
\newtheorem{lem}[thm]{Lemma}
\newtheorem{prop}[thm]{Proposition}
\newtheorem{ax}{Axiom}
\newtheorem{conj}[thm]{Conjecture}

\theoremstyle{definition}
\newtheorem{defn}{Definition}[section]

\theoremstyle{remark}
\newtheorem{rem}{\rm\bfseries{Remark}}[section]
\newtheorem*{notation}{Notation}

\newtheorem{ques}{\rm\bfseries{Question}}[section]
\newtheorem{cons}[rem]{\rm\bfseries{Construction}}
\newtheorem{exm}[rem]{\rm\bfseries{Example}}

%\numberwithin{equation}{section}

      % Includes the theorem environments

\def\OG{\ifmmode{\widetilde{\cal M}_4}\else{$\widetilde{\cal M}_4$}\fi} 
\def\O{\ifmmode{{\mathcal O}}\else{$\mathcal{O}$}\fi} 
\def\H{\ifmmode{{\mathrm H}}\else{${\mathrm H}$}\fi} 
\def\h{\ifmmode{{\mathrm h}}\else{${\mathrm h}$}\fi}
\def\J{\ifmmode{{\mathbb J}}\else{${\mathbb J}$}\fi} 
\def\S{\ifmmode{{\mathbb S}}\else{${\mathbb S}$}\fi} 
\def\K{\ifmmode{{\mathcal K}}\else{${\mathcal K}$}\fi} 
\def\L{\ifmmode{{\mathcal L}}\else{${\mathcal L}$}\fi} 
\def\U{\ifmmode{{\mathcal U}}\else{${\mathcal U}$}\fi} 
\def\F{\ifmmode{{\mathcal F}}\else{${\mathcal F}$}\fi} 
\def\B{\ifmmode{{\mathcal B}}\else{${\mathcal B}$}\fi} 
\def\E{\ifmmode{{\mathcal E}}\else{${\mathcal E}$}\fi} 
\def\M{\ifmmode{{\mathcal M}}\else{${\mathcal M}$}\fi} 
\def\N{\ifmmode{{\mathcal N}}\else{${\mathcal N}$}\fi} 
\def\Br{\ifmmode{{\mathrm{Br}}}\else{${\mathrm{Br}}$}\fi} 
\def\Shah{\ifmmode{\amalg\hspace*{-3.5pt}\amalg}\else{$\amalg\hspace*{-3.5pt}\amalg$}\fi}

\title{Abelian fibred holomorphic symplectic manifolds}
\author[SAWON]{Justin Sawon
%% \\[.2cm] \small{ Dedicated to somebody }
}

%\thanks{Fully supported by T\"UB\.ITAK }

% We write addresses of each author in sequence. These will appear at the
% end of our manuscript.
\address{Department of Mathematics, SUNY at Stony Brook, NY
11794-3651, USA}
\email{sawon@math.sunysb.edu}

\begin{abstract}
We study holomorphic symplectic manifolds which are fibred by abelian
varieties. This structure is a higher dimensional analogue of an
elliptic fibration on a K3 surface. We investigate when a holomorphic
symplectic manifold is fibred in this way, and are led to several
natural conjectures. We then study the geometry of these
fibrations. The expectation is that this point of view will prove
useful in understanding holomorphic symplectic manifolds, and possibly
lead to a classification.
\end{abstract}

%==========================================
%% Do not edit the following command
\volume{9}
%==========================================

\maketitle

\section{Introduction}

Irreducible holomorphic symplectic manifolds are higher dimensional
generalizations of K3 surfaces. It has been roughly twenty years since
Fujiki~\cite{fujiki83} found the first example, the Hilbert scheme
$S^{[2]}$ of two points on a K3 surface $S$, and
Beauville~\cite{beauville83} generalized this to construct two
families, $S^{[n]}$ and the generalized Kummer varieties $K_n$. Since
then there have been other constructions, but only the two examples of
O'Grady~\cite{ogrady99,ogrady00} have given us manifolds which are not
deformations of Beauville's examples. The purpose of this article is
to describe a framework for understanding irreducible holomorphic
symplectic manifolds, which hopefully will lead towards some kind of
classification. The main results are only conjectural, but we will
describe the evidence and motivation behind them, while surveying
special cases which have already been proved.

In studying the moduli space of K3 surfaces, one typically looks at
Kummer surfaces or quartics in $\P^3$, as these are dense but also
relatively easy to understand. However, the structure that will
generalize to higher dimensions is a fibration by abelian
varieties. This suggests that we first review elliptic K3s, which also
happen to be dense in the moduli space. We divide this into three main
steps. Firstly, we need to know which K3 surfaces are
elliptic. Secondly, we describe the family of elliptic K3s which admit
a section. Thirdly, we describe the relation between elliptic K3s
which don't admit sections and their {\em relative Jacobians\/}, which
do. There is nothing new here: for the first step we recall a theorem
of Pjatecki{\u{\i}}-{\v S}apiro and {\v S}afarevi{\v c}~\cite{ps71}
from the 70s (c.f.\ also Section 5 of Kodaira~\cite{kodaira64}), while 
the second and third steps have been well understood for arbitrary
elliptic surfaces for a long time.

Elliptic K3s have base $\mathbb P^1$. There is evidence to suggest
that if the $2n$-dimensional irreducible holomorphic symplectic
manifold $X$ admits a non-trivial fibration, then the fibres must be
abelian varieties and the base must be $\P^n$ (a large part of this
has been proved by Matsushita~\cite{matsushita99}). So the `right'
generalization of an elliptic fibration on a K3 surface appears to be
a fibration by $n$-dimensional abelian varieties over $\P^n$, which we
shall call an {\em abelian fibration\/}. The aim of this paper is to
formulate analogues for abelian fibrations of the three main steps
mentioned above for elliptic K3s.

Firstly, when does an irreducible holomorphic symplectic manifold $X$ 
admit an abelian fibration? Since we are largely interested in
classification up to deformation equivalence, we'd like to know
whether $X$ can always be deformed to have such a structure. Secondly,
can we describe the family of all abelian fibred $X$ which admit
sections? By trying to associate a holomorphic symplectic surface
(ie.\ a K3 surface or complex tori) to such a fibration, we can relate
or possibly even deform it to Beauville's examples. Thirdly, can we
relate abelian fibrations $X$ which don't admit sections to ones that
do? Once again, the goal is to deform $X$ to a more familiar manifold,
which is more amenable to classification.

The programme we will describe tends to suggest that all irreducible
holomorphic symplectic manifolds can be related to Beauville's
examples. By ``related'' we don't simply mean deformed, as O'Grady's
examples would contradict that. Indeed it was through trying to
understand O'Grady's ten-dimensional example and its relation to the
Hilbert scheme $S^{[5]}$ that the author was lead to the ideas in this
paper. The author still feels that a proper understanding of this
relationship will reveal many of the mysteries of holomorphic
symplectic manifolds.

\section{Elliptic K3 surfaces}

Let $S$ be a K3 surface. The group $\H^2(S,\Z)$ is an even unimodular
lattice, with quadratic form $q$ given by intersection pairing. It has
signature $(3,19)$, so from the classification of quadratic forms it
is isomorphic to $L:=3H\oplus 2(-E_8)$. 

\begin{defn}
\label{period}
The {\em period\/} of $S$ is the lattice $(\H^2(S,\Z),q)$ together
with its weight-two Hodge structure
$$\H^2(S,\Z)\otimes\C=\H^2(S,\C)=\H^{2,0}(S,\C)\oplus\H^{1,1}(S,\C)\oplus\H^{0,2}(S,\C).$$
\end{defn}

The holomorphic symplectic form $\sigma$ on $S$ spans
$\H^{2,0}(S,\C)$, and its complex conjugate $\bar{\sigma}$ spans
$\H^{0,2}(S,\C)$. Then $\H^{1,1}(S,\C)$ is given by the orthogonal
complement to $\C\sigma\oplus\C\bar{\sigma}$ with respect to
$q$. Thus after choosing an identification of $(\H^2(S,\Z),q)$ with
$L$, the period is completely determined by the complex line
$[\sigma]\in\P(L\otimes\C)$. This satisfies $q([\sigma])=0$ and
$q([\sigma]+[\bar{\sigma}])>0$, and hence sits inside a quadric $Q$ in 
$\P(L\otimes\C)$. The following result, known as the Global Torelli
Theorem, is originally due to Pjatecki{\u{\i}}-{\v S}apiro and {\v
S}afarevi{\v c}~\cite{ps71} in the algebraic case; the general case is
due to many authors, including Burns, Rapoport, Looijenga, Peters, and
Friedman (for example, see~\cite{lp81} or~\cite{bpv84}).

\begin{thm}
\label{torelli}
If two K3 surfaces $S$ and $S^{\prime}$ have isomorphic periods
$$\psi
:\H^2(S,\Z)\stackrel{\cong}{\longrightarrow}\H^2(S^{\prime},\Z)$$
and moreover the image of the K{\"a}hler cone of $S$ intersects
non-trivially the K{\"a}hler cone of $S^{\prime}$,
$\psi(\K_S)\cap\K_{S^{\prime}}\neq\{0\}$, then there exists a unique
isomorphism
$$f:S\stackrel{\cong}{\longrightarrow}S^{\prime}$$
which induces $\psi$.
\end{thm}

Recall that the K{\"a}hler cone is the open cone in $\H^{1,1}(S,\R)$
of all classes that can be represented by K{\"a}hler forms. If
$\psi(\K_S)\cap\K_{S^{\prime}}=\{0\}$ then $S$ and $S^{\prime}$
will still be isomorphic, but the isomorphism $f$ won't induce
$\psi$. In this case, we first have to compose $\psi$ with a
period-preserving automorphism of the lattice, so that the assumption
on the K{\"a}hler cones is satisfied.

We can define a {\em marked K3 surface\/} to be $S$ together with an
identification of the lattice $(\H^2(S,\Z),q)$ with $L$ (called a {\em
marking\/}). Two marked K3 surfaces are equivalent if there is an
isomorphism between them which commutes with their markings, and we
denote by $\mathcal{T}$ the moduli space of marked K3 surfaces up to
equivalence. One can show that the period map
$$\mathcal{P}:\mathcal{T}\rightarrow Q\subset\P(L\otimes\C),$$
which takes a K3 surface to the complex line $[\sigma]$, is a local
isomorphism, and hence $\mathcal{T}$ is locally isomorphic to a
20-dimensional complex manifold. However, it is not
Hausdorff. Moreover, the actual moduli space $\mathcal{M}$ of K3
surfaces, which is obtained by quotienting $\mathcal{T}$ by the group
of automorphisms of $L$, will be topologically even less
well-behaved.

\subsection{A categorization of elliptic K3 surfaces}
\label{cat_ellK3}

An elliptic K3 surface $S$ will have base $\P^1$. Many kinds of
singular fibres are possible, but generically it will have 24 nodal
$\P^1$s. A fibre $F$ will be a nef divisor with $F^2=0$. If the
fibration has a section $E$, which is necessarily a $(-2)$-curve, then
$(E+F)^2$ also equals zero; however, $(E+F).E=-1$ so $E+F$ is not
nef. Observe that reflection in $E$ recovers the fibre $F$
$$E+F\mapsto (E+F)+((E+F).E)E=F.$$
We recall a result of Pjatecki{\u{\i}}-{\v S}apiro and
{\v S}afarevi{\v c}~\cite{ps71}:

\begin{thm}
\label{ellipticK3}
A projective K3 surface $S$ is elliptic if and only if there exists a
non-trivial divisor $D$ with $D^2=0$.
\end{thm}

\begin{proof}
Here (and throughout) by ``non-trivial'' we mean not linearly equivalent
to the trivial divisor; we write this as $D\neq 0$. Riemann-Roch shows
that either $D$ or $-D$ is effective, so assume that $D$ is
effective. Fix an identification of $(\H^2(S,\Z),q)$ with $L$. The
group of automorphisms of $L$ which fix the period is generated by
\begin{enumerate}
\item the subgroup which preserves the set of nef divisors,
\item the subgroup generated by reflections in $(-2)$-curves,
\item the subgroup generated by $-\mathrm{Id}$.
\end{enumerate}
We will show there is an automorphism belonging to the second of these
subgroups which takes $D$ to a nef divisor. Since $S$ is projective,
let $H$ be an ample divisor and consider $D.H\in\Z$; it is positive as
$D$ is effective. If $D$ is not nef, there exists an irreducible
(effective) divisor $A$ such that $D.A<0$. This can only happen if $A$
is a component of $D$ and $A^2=-2$. Thus $A$ is a $(-2)$-curve and we
can define
$$D_1:=D+(D.A)A$$
to be the reflection in $A$. This reflection preserves the positive
half-cone containing $H$, and therefore $D_1$ is effective. If $D_1$
is not nef, we repeat the procedure. Since
$$D_1.H=D.H+(D.A)(A.H)<D.H$$
we get a decreasing sequence of positive integers $D.H$, $D_1.H$,
$D_2.H$, $\ldots$ and hence the procedure must eventually
terminate. So we can assume that $D$ is nef.

The next step is to show that the linear system $|D|$ contains a
divisor with only one component (possibly non-reduced). This is
achieved by showing that there is at least one divisor in the linear
system which contains an irreducible component $C$ with $C^2=0$; the
ample divisor $H$ is used again at this stage. Then one shows that
$mC\in|D|$ for some $m>0$ (see~\cite{ps71} for details).

Since $C^2=0$, Riemann-Roch shows that $C$ moves in a pencil and we
obtain a rational map $f:S\rightarrow\P^1$; there are no base-points as
$C^2=0$, and hence $f$ is a fibration. Moreover by Bertini's theorem
the generic element of $|C|$ is a smooth irreducible curve, which must
therefore be elliptic.
\end{proof}

\begin{rem}
Earlier Kodaira proved that a K3 surface with Picard group generated
by a non-trivial divisor $D$ with $D^2=0$ must be elliptic (this is
Theorem 15 in~\cite{kodaira64}). Such a surface is necessarily
non-projective, as a projective elliptic K3 surface will contain a
fibre and an ample divisor, and hence will have Picard number at least
two.

Conversely, suppose a K3 surface contains a non-trivial divisor $D$
with $D^2=0$, and has Picard number greater than one. Then the Picard
group, which is a sublattice of $L$, must be indefinite and hence
contains a divisor $E$ with $E^2>0$. The surface is therefore
projective. 

It follows that Kodaira's result is a proof of
Theorem~\ref{ellipticK3} in the non-projective case. Another approach
is to use the argument suggested in Exercise 6 on page 111 of
Beauville's book~\cite{beauville96}. Although this argument is for
projective surfaces, it applies more generally. For instance, we still
have Zariski decomposition for non-projective K3 surfaces by a recent
result of Boucksom~\cite{boucksom02}. The author is grateful to Daniel
Huybrechts for explaining this second proof of
Theorem~\ref{ellipticK3} in the non-projective setting, and to the
referee for directing him to the results in Kodaira's
paper~\cite{kodaira64}.
\end{rem}

\begin{rem}
Suppose we are given $x\in L\backslash\{0\}$ with $x^2=0$, and let $S$
be a K3 surface with a marking
$$\phi:\H^2(S,\Z)\stackrel{\cong}{\longrightarrow}L.$$
Recall the period of $S$ is determined by $[\phi_{\C}(\sigma)]\in
Q\subset\P(L\otimes\C)$; then $\phi^{-1}(x)$ represents a divisor in
$S$ if and only if it has type $(1,1)$, if and only if
$x.\phi_{\C}(\sigma)=0$. For fixed $x$, this gives a hypersurface in
the moduli space $\mathcal{M}$ of all K3 surfaces, so elliptic K3s
occur in 19-dimensional families. Moreover, considering all such
$x\in L$, one can show that the collection of all these families is
dense in $\mathcal{M}$ (the argument is similar to those appearing in
Section 4.7 of~\cite{huybrechts02}). Indeed the denseness of elliptic
K3s in the moduli space of all K3 surfaces was first proved by
Kodaira~\cite{kodaira64}, by deforming to non-projective elliptic K3
surfaces. We use essentially his argument in our proof of
Proposition~\ref{deform_K3} below. 
\end{rem}

In the proof of Theorem~\ref{ellipticK3}, the key point in the first
paragraph is that if there is a curve $A$ on which $D$ is negative,
then $A$ induces an automorphism of the lattice $L$ (preserving the
period) which gives a new divisor $D_1$ with $D_1.A=-D.A$. Moreover,
the process terminates after a finite number of steps. This argument
will fail in higher dimensions. For this reason, let us state and
prove a weaker result that can be more easily generalized (this is
Theorem 18 in~\cite{kodaira64}).

\begin{prop}
\label{deform_K3}
Any K3 surface can be deformed to an elliptic K3.
\end{prop}

\begin{proof}
Deform to a non-projective K3 surface $S$ whose Picard group is generated
by a non-trivial primitive divisor $D$ with $D^2=0$ (thus $D$ must be 
irreducible). That such a deformation exists follows from arguments
similar to those in Section 4.7 of~\cite{huybrechts02}, and only uses
the Local Torelli Theorem identifying $\mathcal{T}$ with $Q$
locally. Either $D$ or $-D$ lies on the boundary of the closure of the
positive cone. For a non-projective K3 surface, the positive cone
coincides with the K{\"a}hler cone (this fact has been generalized in
higher dimensions by Huybrechts~\cite{huybrechts99ii}). Thus either
$D$ or $-D$ is nef, and $S$ contains an irreducible nef divisor with
square zero.

Assuming it is $D$ which is nef, we find
$$\h^2(X,\O(D))=\h^0(X,\O(-D))=0$$
and hence Riemann-Roch shows that $\h^0(X,\O(D))\geq 2$ and $D$ moves
in a pencil. Since $D^2=0$ and $D$ is irreducible, the pencil is
base-point free. It is now fairly easy to show that it is actually
an elliptic fibration over $\P^1$ (this fact has also been partially
generalized in higher dimensions by
Matsushita~\cite{matsushita99,matsushita00}).
\end{proof}

\begin{rem}
In higher dimensions the challenge will be to show that a particular
divisor moves.
\end{rem}

\subsection{Elliptic K3 surfaces which admit sections}

Every elliptic curve over $\C$ can be written as a plane cubic, but
we'd like this description to apply to a $\P^1$ family of elliptic
curves. This is possible provided the family has a section. One way to
understand this is to regard the family as a single elliptic curve $E$ 
over the function field $\C(\P^1)$ of the base; if there exists a
$\C(\P^1)$-valued point then $E$ can be written as a cubic in the
projective plane over $\C(\P^1)$. A $\C(\P^1)$-valued point is
precisely a section of the family.

Suppose that $S$ is an elliptic K3 surface with a section. We assume
that $S$ has generic singular fibres, meaning that the only singular
fibres are nodal or cuspidal $\P^1$s (type $I_1$ or type $II$,
respectively, according to Kodaira's notation~\cite{kodaira63}). If
this were not the case, we should first blow down all the irreducible
components of singular fibres which the section does not pass
through. This will leave us with a singular surface, but at least all
fibres will be irreducible (note that an elliptic K3 surface cannot
have multiple fibres). We can describe $S$ as a family of plane
cubics, where the planes will come from the projectivization of some
rank-three vector bundle $V$ over $\P^1$. By Birkhoff and
Grothendieck's Theorem such a bundle necessarily splits into a direct
sum of line bundles; in fact we find
$$S\subset \P(\O(4)\oplus\O(6)\oplus\O)$$
is given by a cubic equation
$$y^2z=4x^3-axz^2-bz^3$$
where $a$ and $b$ are sections of $\O(8)$ and $\O(12)$
respectively. In this equation $x$, $y$, and $z$ are projections from
$V$ to $\O(4)$, $\O(6)$, and $\O$, so that
$$s:=y^2z-4x^3+axz^2+bz^3$$
is a section of ${\mathrm{Sym}}^3V^*\otimes\O(12)$ whose vanishing cuts
out a cubic curve in each fibre of $\P(V)$. The section of $S$ is
given by $(x,y,z)=(0,1,0)$. Furthermore, there is a $\C^*$-action
$$a\mapsto\la^4a\qquad\mbox{and}\qquad b\mapsto\la^6b$$
such that all pairs $(a,b)$ in the same orbit give isomorphic surfaces
$S$.

The sections $a$ and $b$ belong to 9 and 13-dimensional vector
spaces, respectively. Taking into account the $\C^*$-action and the
3-dimensional family ${\mathrm{PGL}}(2,\C)$ of automorphisms of the
base $\P^1$, we see that the moduli space of elliptic K3 surfaces
which admit a section is 18-dimensional. Moreover, this moduli space
is irreducible by a theorem of Seiler (proved in Friedman and
Morgan~\cite{fm94}, for example).

The description of an elliptic surface as a family of cubics is known
as a {\em Weierstra{\ss} model\/}, and exists for any elliptic surface
which admits a section (see~\cite{bpv84}, \cite{fm94}, \cite{is96}, or
almost any book on complex surfaces).

\subsection{Elliptic K3 surfaces which don't admit sections}
\label{K3_nosection}

Next we wish to consider elliptic K3 surfaces which don't necessarily
admit sections. Once again, we will assume that $S$ has generic
singular fibres (nodal or cuspidal $\P^1$s). Let $i:S_t\hookrightarrow
S$ be the inclusion of a smooth fibre, and let $i_*L$ be the
push-forward of a degree-zero line bundle $L$ on $S_t$. Associated to
$S$ is its {\em relative Jacobian\/} $J$, defined as the moduli space
of pure semi-stable sheaves on $S$ with the same Hilbert polynomial as
$i_*L$. Since $S$ has generic singular fibres (in particular, all
fibres are irreducible), every element of $J$ is the push-forward of a
degree-zero torsion-free sheaf on some fibre. Thus there is a
projection from $J$ to $\P^1$:
$$\begin{array}{ccccc}
S & & & & J \\
 & p_S\searrow & & \swarrow p_J & \\
 & & {\P}^1 & & \\
\end{array}$$
Each fibre $J_t={\mathrm{Jac}}(S_t)$ is isomorphic to the
corresponding $S_t$, and in fact $S$ and $J$ are locally isomorphic as
fibrations. However, the relative Jacobian $J$ has a natural section
$s_0$ given by the flat family $\O_S$ (thus
$s_0(t)=\O_S|_{S_t}=\O_{S_t}\in J_t$) whereas $S$ itself may not have
a section. Let $\{U_i\}$ be an open cover of $\P^1$ with $S_i\cong
J_i$ over $U_i$, where we denote $p_S^{-1}(U_i)$ and $p_J^{-1}(U_i)$
by $S_i$ and $J_i$ respectively. The isomorphism is not unique, but is
given by the choice of a local section of $S$: if $s_i:U_i\rightarrow
S_i$ is a local section, there exists a local isomorphism
$\phi_i:S_i\rightarrow J_i$ taking $s_i$ to $s_0$. If the $s_i$ patch
together to give a global section, then the $\phi_i$ patch together to
give a global isomorphism $S\cong J$. Hence we have

\begin{prop}
The elliptic fibration $S$ is isomorphic to its relative Jacobian $J$ 
if and only if it admits a section.
\end{prop}

In general, we can regard $S$ as being a fibration which is locally
the same as $J$, but stuck together in a different way. In other words
$S$ is a {\em torsor\/} over $J$, and hence is classified by an
element of $\H^1(\P^1,\B)$ where $\B$ is the sheaf of local
holomorphic sections of $p_J:J\rightarrow \P^1$ (or equivalently, the
sheaf of translations in fibres, as a local section $s$ gives a family
of translations by $s-s_0$). To see this explicitly, observe that
above an overlap $U_i\cap U_j$ the local isomorphisms
$\phi_i:S_i\rightarrow J_i$ and $\phi_j:S_j\rightarrow J_j$ differ by
a translation, and hence $\phi_i\circ\phi_j^{-1}$ is a local section
of $\B$. Together these give a {\v C}ech 1-cocycle. 

\begin{defn}
The {\em analytic Tate-Shafarevich group\/} $\Shah^{\mathrm{an}}(J)$
is the group which classifies, up to isomorphism, elliptic fibrations
$S$ with relative Jacobian $J$.
\end{defn}

The group $\Shah^{\mathrm{an}}(J)$ classifies {\em all\/} such $S$,
including those with multiple fibres. Since these additional surfaces
will not be K3 surfaces, we ignore them, and we are left with
$\H^1(\P^1,\B)$. We can use the Leray spectral sequence to show that 
$$\H^1(\P^1,\B)\cong\H^2(J,\O^*).$$
The latter is the {\em analytic Brauer group of $J$\/}; it is
one-dimensional, but with a strange topology. More precisely, the
exponential exact sequence on $J$ gives
$$0\rightarrow\H^1(J,\O^*)\rightarrow\H^2(J,\Z)\rightarrow\H^2(J,\O)\rightarrow\H^2(J,\O^*)\rightarrow
0.$$
The first group is the Picard group ${\mathrm{Pic}}(J)$, the second is
the rank 22 lattice $L$, and the third is isomorphic to $\C$. Thus
$\H^2(J,\O^*)$ is the quotient of $\C$ by a lattice which can vary
from rank 20, when $J$ has Picard number $\rho=2$, to rank 2, when
$\rho=20$ (note that $\rho\geq 2$ since $J$ is elliptic with a section,
and so always has at least two independent divisors, namely a fibre
and a section).

The extra dimension of the Brauer group, added to the 18 of the moduli
space of elliptic K3 surfaces which admit sections, gives a
19-dimensional moduli space of elliptic K3 surfaces, as expected.

\begin{rem}
The Brauer group $\H^2(J,\O^*)$ is the space of {\em gerbes\/} on
$J$. In Section 6 we shall give a description, due to C{\u
a}ld{\u a}raru, of how gerbes arise in the classification of elliptic
fibrations $S$ with given relative Jacobian $J$. 
\end{rem}

\begin{rem}
\label{torsion}
If $S$ is algebraic then an ample divisor induces a multi-valued
section of the fibration $S\rightarrow \P^1$. Suppose $S$ corresponds
to an element $\a$ of the Brauer group $\H^2(J,\O^*)$ of its relative
Jacobian. If $S$ admits a $k$-valued section then one can show that
$\a$ must be $k$-torsion, and conversely. The {\em algebraic Brauer
group\/} consists of the torsion elements of $\H^2(J,\O^*)$
(equivalently, we can define it as $\H_{\mathrm{\acute{e}t}}^2(J,\O^*)$,
where we calculate cohomology in the {\'e}tale topology). Thus
algebraic elliptic fibrations correspond to elements of the algebraic
Brauer group, which is dense in the analytic Brauer group.
\end{rem}

Once again, this theory is not specific to K3 surfaces, but applies to
general elliptic surfaces, and was first developed in this form by
Kodaira~\cite{kodaira63} (alternatively see~\cite{bpv84}, \cite{fm94},
or~\cite{is96}).

This concludes what we want to say about elliptic K3 surfaces. In
summary, the three main points are
\ben
\item any K3 surface can be deformed to an elliptic K3 surface,
\item elliptic K3 surfaces which admit sections are described by an
18-dimensional {\em connected\/} moduli space,
\item an elliptic K3 surface which does not admit a section can be
connected to its relative Jacobian by a one-dimensional family.
\een
In particular, all K3 surfaces are deformation equivalent. In the rest
of this article we will formulate generalizations of these statements
for higher dimensional irreducible holomorphic symplectic manifolds.

\section{Irreducible holomorphic symplectic manifolds}

K3 surfaces are compact, K{\"a}hler, simply-connected, and have
non-degenerate holomorphic two-forms. All of these notions generalize
to higher dimensions.

\begin{defn}
A {\em holomorphic symplectic manifold\/} $X$ is a compact K{\"a}hler
manifold which admits a non-degenerate holomorphic two-form
$\s\in\H^0(X,\La^2 T^*)$, known as the {\em holomorphic symplectic
form\/}. If $X$ is simply-connected and cannot be written as a product
$Y\times Z$ of two manifolds, then we shall say $X$ is {\em
irreducible\/}.
\end{defn}

Suppose $X$ has dimension $2n$. By non-degeneracy, $\s^n$ is a
trivialization of the canonical bundle $\K=\La^{2n}T^*$. Therefore
$c_1(T)=0$ and by Yau's theorem $X$ admits a Ricci-flat metric. In
fact, this metric is {\em hyperk{\"a}hler\/}, and thus compact
hyperk{\"a}hler manifolds are really the differential geometric
equivalents of holomorphic symplectic manifolds.

\subsection{Examples}
\label{examples}

The following two families of examples are due to
Beauville~\cite{beauville83}.

\begin{exm}
\label{HilbK3}
Obviously a K3 surface $S$ is an irreducible holomorphic symplectic
manifold. The {\em Hilbert scheme of $n$ points on $S$\/}, denoted
$S^{[n]}$ is the moduli space parameterizing length $n$
zero-dimensional subschemes of $S$. A typical element consists of $n$
distinct unordered points of $S$, though the behaviour when some of
the points collide is more complicated. For example, $S^{[2]}$ is
given by blowing up the diagonal in $S\times S$ and then quotienting
by the involution which exchanges the two factors (this was the first
higher dimensional holomorphic symplectic manifold, discovered by
Fujiki~\cite{fujiki83}). The Hilbert scheme $S^{[n]}$ is a smooth
resolution of ${\mathrm{Sym}}^nS$. It is an irreducible holomorphic
symplectic manifold of dimension $2n$. 
\end{exm}

\begin{exm}
\label{genKum}
Let $T$ be an abelian surface, or more generally a complex torus of
dimension two. The Hilbert scheme $T^{[n+1]}$ is a holomorphic
symplectic manifold, but it is not irreducible. There is a map $\pi$
to $T$ given by composing the map
$T^{[n+1]}\rightarrow{\mathrm{Sym}}^{n+1}T$ with the group structure
on $T$. The fibres are all isomorphic. Define the {\em generalized
Kummer variety\/} $K_n$ to be $\pi^{-1}(0)$. It is an irreducible
holomorphic symplectic manifold of dimension $2n$. When $n=1$, we find
$K_1$ is a Kummer K3 surface, whence the name.
\end{exm}

The following example is due to Mukai~\cite{mukai84,mukai87}.

\begin{exm}
\label{mukaiK3}
Let $S$ be a K3 surface. We saw that there is a quadratic form $q$ on 
$\H^2(S,\Z)$ given by intersection pairing. We extend this to a
quadratic form on $\H^{\bullet}(S,\Z)$ by defining
$$q(v,w):=\int_S(-v_0w_4+v_2w_2-v_4w_0)$$
where $v_i$ and $w_i\in\H^i(S,\Z)$. For a sheaf $\E$ on $S$, define
the {\em Mukai vector\/} $v(\E)\in\H^{\bullet}(S,\Z)$ to be the
product ${\mathrm{ch}}(\E){\mathrm{Td}}^{1/2}$ of the Chern character
of $\E$ with the square root of the Todd class of $S$. For example, if
$\E$ is locally free of rank $r$, with Chern classes $c_1$ and $c_2$,
then 
$$v(\E)=(r,c_1,r+c_1^2/2-c_2).$$
For $v\in\H^{\bullet}(S,\Z)$, define $\M^{\mathrm{s}}(v)$ to be the moduli
space of stable sheaves $\E$ on $S$ with Mukai vector $v(\E)=v$. Then
$\M^{\mathrm{s}}(v)$ is smooth of dimension $2n=q(v,v)+2$ and has a
holomorphic symplectic form. It is known that for various choices of
$v$, $\M^{\mathrm{s}}(v)$ is birational to the Hilbert scheme
$S^{[n]}$: this is immediate for rank one, was proved for various rank
two moduli spaces in~\cite{gh96}, and for arbitrary rank and elliptic
K3s it was proved in~\cite{ogrady97} (see also Section 11.3
in~\cite{hl97}). However, in general the moduli space of stable
sheaves need not be compact. To compactify we have to add semi-stable
sheaves, and this may result in singularities.

However, if $v$ is primitive then every semi-stable sheaf $\E$ with
Mukai vector $v(\E)=v$ will automatically be stable. So in this case
$\M^{\mathrm{s}}(v)$ is an irreducible holomorphic symplectic
manifold. On the other hand, under these conditions it is also known
that $\M^{\mathrm{s}}(v)$ is a deformation of the Hilbert scheme
$S^{[n]}$, where $2n=q(v,v)+2$. This follows from the work of
G{\"o}ttsche, Huybrechts, and
O'Grady~\cite{gh96,huybrechts99i,ogrady97} (the most general statement
is due to Yoshioka~\cite{yoshioka99}).
\end{exm}

Mukai's results also apply to the moduli space of sheaves on an
abelian surface. Unfortunately, when we can show that the moduli space
is smooth and compact, it also happens to be a deformation of the
generalized Kummer variety (see Yoshioka~\cite{yoshioka01}, for
example).

The first known example not deformation equivalent to a Hilbert scheme
$S^{[n]}$ or a generalized Kummer variety $K_n$ is due to
O'Grady~\cite{ogrady99}.

\begin{exm}
\label{ogrady}
Consider the previous example with Mukai vector $v=(2,0,-2)$. This
vector is not primitive, and hence there exist strictly semi-stable
sheaves which we add to $\M^{\mathrm{s}}(v)$ to get the {\em
singular\/} moduli space of semi-stable sheaves
$\M^{\mathrm{ss}}(v)$. O'Grady proved that $\M^{\mathrm{ss}}(v)$
admits a symplectic desingularization $\widetilde{\M}$, which is a
ten-dimensional irreducible symplectic manifold. He also showed that
$b_2(\widetilde{\M})\geq 24$, and hence $\widetilde{\M}$ is not a
deformation of the Hilbert scheme $S^{[5]}$ or the generalized Kummer
variety $K_5$, which have second Betti numbers equal to 23 and 7,
respectively.
\end{exm}

The only other known example comes from a similar construction
starting with an abelian surface (see O'Grady~\cite{ogrady00} for
details). It is six-dimensional and has $b_2=8$.

\subsection{Moduli spaces}

Moduli spaces of irreducible holomorphic symplectic manifolds have
been extensively studied~\cite{huybrechts99i,huybrechts02}. Although
the Global Torelli Theorem is false in general, the theory is
nonetheless very similar to the K3 case. For an arbitrary irreducible
holomorphic symplectic manifold $X$, there exists a quadratic form
$q_X$ on $\H^2(X,\Z)$ known as the Beauville-Bogomolov
form~\cite{beauville83}; this generalizes the intersection pairing on
$\H^2(S,\Z)$ when $S$ is a K3 surface. So $(\H^2(X,\Z),q_X)$ is once
again a lattice that we write abstractly as $(\Gamma,q_{\Gamma})$ (or
simply as $\Gamma$). As in Definition~\ref{period}, we can define the
period of $X$ as $(\H^2(X,\Z),q_X)$ together with its weight-two Hodge
structure, and after identifying $\H^2(X,\Z)$ with $\Gamma$, this is
once again determined by the complex span of the holomorphic
symplectic form $[\sigma]\in\P(\Gamma\otimes\C)$. 

As with K3 surfaces, there is a period map 
$$\mathcal{P}_{\Gamma}:\mathcal{T}_{\Gamma}\rightarrow
Q_{\Gamma}:=\{x\in\P(\Gamma\otimes\C)|q_{\Gamma}(x)=0,\mbox{
}q_{\Gamma}(x+\bar{x})>0\}$$
from the moduli space $\mathcal{T}_{\Gamma}$ of {\em marked\/}
irreducible holomorphic symplectic manifolds to a quadric in
$\P(\Gamma\otimes\C)$. Beauville showed the period map is a local
isomorphism, and hence $\mathcal{T}_{\Gamma}$ is locally isomorphic to
a complex manifold of dimension $b_2-2$ (where $b_2$ is the second
Betti number of $X$), though as with K3 surfaces
$\mathcal{T}_{\Gamma}$ is not Hausdorff. 

Huybrechts~\cite{huybrechts99i} proved surjectivity of the
period map, and also showed that non-separated points in
$\mathcal{T}_{\Gamma}$ correspond to birational (marked) manifolds $X$ 
and $X^{\prime}$. After quotienting $\mathcal{T}_{\Gamma}$ by
automorphisms of $\Gamma$ to obtain $\mathcal{M}_{\Gamma}$ the
converse is also true: (unmarked) $X$ and $X^{\prime}$ are birational
if and only if they correspond to non-separated points in
$\mathcal{M}_{\Gamma}$. This also implies the important result that
two birational holomorphic symplectic manifolds $X$ and $X^{\prime}$
are deformation equivalent.

Finally, let us also note that the Global Torelli Theorem as
formulated for K3 surfaces (Theorem~\ref{torelli}) is false in higher
dimensions: Namikawa~\cite{namikawa01} has constructed a
counter-example consisting of two non-birational holomorphic
symplectic manifolds with isomorphic periods.

\subsection{Abelian fibrations}

\begin{defn}
By {\em abelian fibration\/} on a $2n$-dimensional irreducible
holomorphic symplectic manifold $X$ we shall mean the structure of a
fibration over $\P^n$ whose generic fibre is a smooth abelian variety
of dimension $n$.
\end{defn}

This shall be our higher dimensional analogue of elliptic
fibrations on K3 surfaces. At first sight, the definition may appear
to be unnecessarily restrictive. For example, maybe we should allow
the base to be a more general $n$-fold than $\P^n$, or even to have 
dimension different to $n$. However, this is more-or-less the only
fibration structure that can exist on an irreducible holomorphic
symplectic manifold, by the following result of
Matsushita~\cite{matsushita99,matsushita00}.

\begin{thm}
\label{matsushita}
For projective $X$, let $f:X\rightarrow B$ be a proper surjective
morphism such that the generic fibre $F$ is connected. Assume that $B$
is smooth and $0<{\mathrm{dim}}B<{\mathrm{dim}}X$. Then
\ben
\item $F$ is an abelian variety up to a finite unramified cover,
\item $B$ is $n$-dimensional and has the same Hodge numbers as $\P^n$, 
\item the fibration is Lagrangian with respect to the holomorphic
symplectic form.
\een
\end{thm}

In particular, if $X$ is 4-dimensional, we can use the
Castelnuovo-Enriques classification of surfaces to deduce that the
generic fibre is an abelian surface and the base is $\P^2$
(this was also proved by Markushevich in~\cite{markushevich95}, but
under the assumption that the fibration is Lagrangian). Assuming only
that $B$ is normal, Matsushita also deduced slightly weaker results.

Matsushita's theorem provides some justification for assuming the base
must be $\P^n$. Moreover, most of (perhaps all of) the examples of
irreducible holomorphic symplectic manifolds described in
Section~\ref{examples} can be deformed to abelian fibrations.

\begin{exm}
Both Examples~\ref{HilbK3} and~\ref{genKum}, the Hilbert scheme of
points on a K3 surface and the generalized Kummer variety, are abelian
fibrations when the underlying K3 surface $S$ or complex tori $T$,
respectively, is an elliptic surface. For example, if
$f:S\rightarrow\P^1$ is the fibration on $S$, we get an induced
fibration 
$$f^{[n]}:S^{[n]}\rightarrow{\mathrm{Sym}}^n\P^1\cong\P^n$$
on $S^{[n]}$. However, the fibres in this case are products of $n$
elliptic curves: still $n$-dimensional abelian varieties, but quite a
special case. A similar thing happens for the generalized Kummer
variety. We regard these abelian fibrations as rather degenerate;
there are `better' abelian fibrations on these manifolds, as the next
two examples show.
\end{exm}

\begin{exm}
\label{abel_fibr_K3}
Let $S$ be a K3 surface which contains a smooth genus $g$ curve
$C$. Let $L$ be a degree $g$ line bundle on a smooth curve $D$ in
the linear system $|C|$, and consider $i_*L$ where $i$ is the
inclusion of $D$ in $S$. Let $Z$ be the moduli space of rank-one
torsion sheaves on $S$ which have the same type as $i_*L$. In other
words, $Z$ is the Mukai moduli space $\M^{\mathrm{s}}(0,[C],1)$. Since 
$v=(0,[C],1)$ is primitive, the moduli space is smooth and
compact. Clearly it has the structure of an abelian fibration, with
base given by the linear system $|C|\cong\P^g$ and fibre over $D\in
|C|$ given by ${\mathrm{Pic}}^gD$:
$$\begin{array}{ccc}
	{\mathrm{Pic}}^g & \hookrightarrow & Z \\
                         &                 & \downarrow \\
                         &                 & |C|\cong\P^g \\
  \end{array}$$
From the comments in Example~\ref{mukaiK3}, we know that $Z$ is
deformation equivalent to the Hilbert scheme $S^{[g]}$. It is also
birational to the Hilbert scheme. To see this, note that by
Riemann-Roch a generic degree $g$ line bundle on a genus $g$ curve has
a unique section up to scale, which vanishes at precisely $g$
points. This gives a map from a (Zariski) open subset $U$ of $Z$ to
$S^{[g]}$. To get a map in the opposite direction, recall that for
$g\geq 3$, $S$ is embedded in the dual $(\P^g)^{\vee}$ of the linear
system $|C|$ (for $g=2$, $S$ is a double cover of $(\P^2)^{\vee}$
branched over a sextic). A generic collection of $g$ points on $S$
defines a hyperplane, which intersects $S$ in a curve $D\in |C|$. The
$g$ points lie on this curve and define a degree $g$ line bundle on
$D$. Hence we have a map from an open subset of $S^{[g]}$ to $Z$.
\end{exm}

\begin{rem}
The birational map in Example~\ref{abel_fibr_K3} is a {\em generalized
Mukai flop\/}, and was studied by Markman~\cite{markman01}. Also,
recall that two birational holomorphic symplectic manifolds are
automatically deformation equivalent, by the results of
Huybrechts~\cite{huybrechts99i}.
\end{rem}

\begin{exm}
\label{abel_fibr_genK}
For generalized Kummer varieties we can do something similar. We begin
with an abelian surface $T$ which contains a genus $g+2$ curve
$C$. Let $Y^{\prime}$ be the moduli space of rank-one torsion sheaves
on $T$, which have the same type as push-forwards of degree $g+2$ line
bundles on smooth curves in the linear system $|C|$. This is the
Mukai moduli space $\M^{\mathrm{s}}(0,[C],1)$ of sheaves on $T$, and
it is smooth and compact. However, $Y^{\prime}$ is not
irreducible. Since a degree $g+2$ line bundle on a genus $g+2$ curve
$C$ determines, generically, $g+2$ points on $C$, we obtain a map
$\pi$
$$Y^{\prime}\rightarrow {\mathrm{Sym}}^{g+2}T\rightarrow T$$
where the second map is given by the group structure on $T$. This is
just the Albanese map of $Y^{\prime}$, and the fibres are all
isomorphic. Then $Y:=\pi^{-1}(0)$ is an irreducible holomorphic
symplectic manifold of dimension $2g$. Moreover, it is a deformation
of the generalized Kummer variety $K_g$, by results of
Yoshioka~\cite{yoshioka01}.

We will show that $Y$ is an abelian fibration. Since the linear system
of curves $|C|$ on $T$ is $g$-dimensional, $Y^{\prime}$ has a
fibration structure
$$\begin{array}{ccc}
	{\mathrm{Pic}}^{g+2} & \hookrightarrow & Y^{\prime} \\
                         &                 & \downarrow \\
                         &                 & |C|\cong\P^g \\
  \end{array}$$
which induces a fibration of $Y$
$$\begin{array}{ccc}
	Y\cap{\mathrm{Pic}}^{g+2} & \hookrightarrow & Y \\
                         &                 & \downarrow \\
                         &                 & |C|\cong\P^g \\
  \end{array}$$
Moreover, when we restrict $\pi$ to the fibre ${\mathrm{Pic}}^{g+2}D$
of $Y^{\prime}\rightarrow\P^g$, $D\in |C|$, it still surjects onto
$T$. Therefore the induced fibre of $Y$
$$Y\cap{\mathrm{Pic}}^{g+2}D=\mathrm{ker}(\pi|_{{\mathrm{Pic}}^{g+2}D}:{\mathrm{Pic}}^{g+2}D\rightarrow T)$$
will be a $g$-dimensional abelian variety. Therefore the generalized
Kummer variety $K_g$ can be deformed to an abelian fibration $Y$.
\end{exm}

\begin{rem}
\label{polarizations}
In Example~\ref{abel_fibr_K3} the fibres of $Z$ (which is a
deformation of $S^{[g]}$) are principally polarized abelian varieties,
as they are Jacobians of curves. In Example~\ref{abel_fibr_genK} the
fibres of $Y$ (which is a deformation of $K_g$) are not Jacobians;
they are given by
$$Y_D:=\mathrm{ker}(\pi|_{{\mathrm{Pic}}^{g+2}D}:{\mathrm{Pic}}^{g+2}D\rightarrow
T)$$
where $D\in |C|$. Assume that $T$ is polarized by the curve $C$, and
let the type of this polarization be $(d_1,d_2)$ where $d_1|d_2$ (see
Birkenhake and Lange~\cite{bl92}). If we assume that $C$ is reduced,
then from
$$d_1d_2=\frac{C.C}{2}=g+1.$$
we conclude that $(d_1,d_2)=(1,g+1)$. Therefore the (complementary)
polarization of $Y_D$ has type $(1,\ldots,1,g+1)$. In particular, the
fibres of $Y$ are never principally polarized. We shall return to this
point later.
\end{rem}

\begin{exm}
\label{abel_fibr_ogrady}
Let $S$ be a K3 surface which is the double cover of $\P^2$ branched
over a sextic. The pull-back of a conic from $\P^2$ is a genus five
curve $C$ on $S$. In~\cite{ogrady99}, O'Grady showed that an open
subset of his moduli space $\widetilde{\M}$ (from
Example~\ref{ogrady}) is birational to an open subset of the Mukai
moduli space $\M^{\mathrm{s}}(0,[C],2)$ of sheaves on $S$. The latter
is the moduli space of rank-one torsion sheaves with the same type as 
push-forwards of degree $6$ line bundles on smooth curves in $|C|$.
In this case $v=(0,[C],2)$ is not primitive, as curves in $|C|$ can be
non-reduced (for example, the pull-back of a double line). So a priori,
it is not clear that the moduli space $\M^{\mathrm{s}}(0,[C],2)$ can
be completed to a smooth irreducible holomorphic symplectic
manifold. More precisely, adding strictly semi-stable sheaves will
give a singular space $\M^{\mathrm{ss}}(0,[C],2)$, and it is not known 
whether a symplectic desingularization of this space exists (though
presumably one can use the same kind of desingularization as for
$\widetilde{\M}$).

In any case, the incomplete space $\M^{\mathrm{s}}(0,[C],2)$ obviously
has the structure of an abelian fibration 
$$\begin{array}{ccc}
	{\mathrm{Pic}}^6 & \hookrightarrow & {\M}^{\mathrm{s}}(0,[C],2) \\
                         &                 & \downarrow \\
                         &                 & |C|\cong\P^5 \\
  \end{array}$$
Since birational varieties are deformation equivalent, it appears
probable that O'Grady's space $\widetilde{\M}$ can be deformed to an 
abelian fibration. 
\end{exm}

There is likely to be a similar description for O'Grady's other
space~\cite{ogrady00}. If we could make these arguments rigorous, it
would imply that all the {\em currently known\/} examples of
irreducible holomorphic symplectic manifolds can be deformed to
abelian fibrations. This is the first step in a (conjectural)
three-part programme to understand irreducible holomorphic symplectic
manifolds, which we now summarize
\ben
\item given $X$, deform it to an abelian fibration,
\item classify those abelian fibrations which admit sections,
\item connect abelian fibrations which don't admit sections to
corresponding fibrations which do by one-dimensional families.
\een
We shall address each of these problems in the next three sections.

\section{Deforming to abelian fibrations}

In this section, we address the question of whether an arbitrary
irreducible holomorphic symplectic manifold can be deformed to an
abelian fibration.

\subsection{Divisors on abelian fibrations}

By Theorem~\ref{ellipticK3} a (projective) K3 surface $S$ is elliptic
if and only if it admits a non-trivial divisor whose square is
zero. A fibre $F$ of the elliptic fibration is an example of such a
divisor, though as we saw at the beginning of Section~\ref{cat_ellK3},
not the only example. What special divisors exist on an abelian fibred
irreducible holomorphic symplectic manifold $X^{2n}$? The most obvious
choice is the pull-back of a hyperplane from the base $\P^n$; call 
this $D$. In fact, in the most generic situation the Picard group of
$X$ will be generated by $D$. 

Observe that $D^n$ is non-trivial, as it is rationally equivalent to a
fibre of $X$, but $D^{n+1}=0$. For a divisor $E$, the largest integer
$k$ such that $E^k$ is non-trivial is known as the {\em numerical
Iitaka dimension\/} of $E$ (this is sometimes referred to as the {\em
numerical Kodaira dimension\/}, but following Esnault and
Viehweg~\cite{ev92}, we reserve such terminology for the canonical
divisor). Thus $D$ has numerical Iitaka dimension $n$.

For an arbitrary irreducible holomorphic symplectic manifold $X^{2n}$,
Verbitsky~\cite{verbitsky96} proved that the symmetric product
${\mathrm{Sym}}^{\bullet}\H^2(X,\C)$ of the second cohomology group
injects into the cohomology ring $\H^{\bullet}(X,\C)$ up to the middle
dimension $\H^{2n}(X,\C)$. It follows that given an arbitrary
non-trivial divisor $D$ on $X$, $D^n$ must be non-trivial. On the
other hand, Verbitsky also showed that $D^{n+1}=0$ if and only if
$q_X(D)=0$, where $q_X$ is the Beauville-Bogomolov quadratic
form. Fujiki's formula (see~\cite{huybrechts99i}, for example) says
that
\begin{eqnarray}
\label{fujiki}
\int_XD^{2n} & = & c_Xq_X(D)^n
\end{eqnarray}
for all divisors, where $c_X$ is a positive real scalar
(independent of $D$) known as the {\em Fujiki constant\/}. So if
$D^k=0$ for $k\geq n+1$, then $q_X(D)=0$ and hence $D^{n+1}=0$. In
summary, a divisor $D$ on $X^{2n}$ will have numerical Iitaka
dimension equal to
\ben
\item 0 if $D$ is trivial,
\item $n$ if $q_X(D)=0$,
\item $2n$ otherwise.
\een

We saw above that an abelian fibration $X$ has a non-trivial divisor
whose square is zero with respect to $q_X$, just as an elliptic K3 surface
$S$ has a divisor whose square is zero with respect to the
intersection pairing. Could the converse be true?

\begin{exm}
\label{nef_is_nec}
Let $g:S\rightarrow\P^2$ be a K3 surface which is a double cover of
the plane branched over a sextic. The pull-back of a line from $\P^2$
is a genus two curve $C$, and we saw in Example~\ref{abel_fibr_K3}
that the moduli space $\M^{\mathrm{s}}(0,[C],1)$ is birational to the
Hilbert scheme $S^{[2]}$. The birational transformation is an example
of a Mukai flop~\cite{mukai84}. It is given by first blowing-up the
locus
$$\{\K_D|D\in|C|\}\cong(\P^2)^{\vee}\subset\M^{\mathrm{s}}(0,[C],1)$$
consisting of canonical bundles (note that $|C|$ is the set of lines
in the original $\P^2$, and hence this locus is isomorphic to the dual
plane $(\P^2)^{\vee}$). The exceptional locus of the blow-up is
isomorphic to $$\{(x,l)\in\P^2\times(\P^2)^{\vee}|x\in l\}$$
and can be blown-down in a different direction to produce the locus
$$G:=\{g^{-1}(x)|x\in\P^2\}\cong\P^2\subset S^{[2]}$$
(more precisely, $G$ is the closure of the set of $g^{-1}(x)$ for $x$
in the complement of the sextic branching curve).

Now $\M^{\mathrm{s}}(0,[C],1)$ is an abelian fibration over
$(\P^2)^{\vee}$. Choose a point $w\in\P^2$; it will corresponds to a
hyperplane in $(\P^2)^{\vee}$ and thus we obtain a divisor $D$ in
$\M^{\mathrm{s}}(0,[C],1)$ whose square is zero. Since birational maps
between holomorphic symplectic manifolds are isomorphisms in
codimension two (see for instance~\cite{huybrechts99i} or Huybrechts'
notes in~\cite{ghj02}), there must be a corresponding divisor
$D^{\prime}$ in $S^{[2]}$ and it is easy to see that
$$D^{\prime}=\{\xi\in S^{[2]}|\mbox{$w$, $y$, and $z$ are collinear,
where }\{y,z\}=g(\mathrm{supp(\xi)})\}$$
(once again, $D^{\prime}$ should really be defined as the closure of
the set of such $\xi$ with $y$ and $z$ distinct). From this
description, one can see that the rational equivalence class of
$(D^{\prime})^2$ is not base-point free; indeed $D^{\prime}$ always
contains the locus
$$G=\{g^{-1}(x)|x\in\P^2\}\cong\P^2$$
regardless of the choice of $w$. So $(D^{\prime})^2$ does not
represent a fibre and the structure of an abelian fibration is lost
under the birational transform (of course, we have only shown that
$D^{\prime}$ is not the pull-back of a hyperplane from the base of an
abelian fibration; $S^{[2]}$ might still have a fibred structure
unrelated to $D^{\prime}$).

Why does $D^{\prime}$ not induce an abelian fibration on $S^{[2]}$?
This appears to be because $D^{\prime}$ is not nef: any curve in the
base-point locus $G$ will intersect $D^{\prime}$ negatively. On the
other hand, $D$ is nef: a curve in
$\M^{\mathrm{s}}(0,[C],1)\rightarrow(\P^2)^{\vee}$ is either contained
in a fibre and has zero intersection with $D$, or is a finite cover of
a curve in $(\P^2)^{\vee}$ and intersects $D$ positively.
\end{exm}

This example highlights two important points. Firstly, if $X$ contains
a divisor $D$ with $q_X(D)=0$, and we want to show that $X$ is an
abelian fibration, then it will probably help to assume that $D$ is
nef. Secondly, if $q_X(D)=0$ but $D$ is not nef, then it may still be
possible to find a birational model $X^{\prime}$ of $X$ which contains
a nef divisor whose square is zero, and which is an abelian
fibration. We state this in the following conjectures.

\begin{conj}
\label{nef_sqzero}
An irreducible holomorphic symplectic manifold $X$ is an abelian
fibration if and only if it contains a non-trivial nef divisor $D$
whose square with respect to the Beauville-Bogomolov quadratic form is
zero, $q_X(D)=0$.
\end{conj}

\begin{rem}
Clearly an abelian fibration contains such a divisor, namely the
pull-back of a hyperplane from the base. We will try to argue (in this
subsection and the next) that the right kind of vanishing theorems
will imply the converse.
\end{rem}

\begin{conj}
\label{just_sqzero}
Any irreducible holomorphic symplectic manifold $X^{\prime}$ which
contains a non-trivial divisor $D^{\prime}$ with $q_X(D^{\prime})=0$
is birational to an abelian fibration $X$. Moreover, if $D$ is the
divisor corresponding to $D^{\prime}$, then there exists a
period-preserving automorphism $\tau$ of the lattice $\H^2(X,\Z)$ such
that $\tau(D)$ is nef.
\end{conj}

\begin{rem}
In Example~\ref{nef_is_nec} the divisor $D$ corresponding to
$D^{\prime}$ is already nef, but as we saw in the proof of
Theorem~\ref{ellipticK3} for elliptic K3s, the automorphism $\tau$ is 
necessary in general.
\end{rem}

Conjecture~\ref{just_sqzero} is based on observations like in
Example~\ref{nef_is_nec} and on our insight from the K3 case. Let us
look instead at Conjecture~\ref{nef_sqzero}. We saw already that a
non-trivial divisor $D$ with $q_X(D)=0$ has numerical Iitaka dimension
$n$. Consider the map to projective space
$$\phi_{mD}:X\rightarrow\P(\H^0(X,\O(mD)))$$
given by the sections of $\O(mD)$, for $m>0$. The dimension of the
image of $\phi_{mD}$, in the limit as $m\rightarrow\infty$, is known
as the {\em Iitaka dimension\/} of $D$. We want to show that the
Iitaka dimension of $D$ is also $n$, as it is the map $\phi_{mD}$
which ought to give the abelian fibration with base space $\P^n$ (of
course we also have to prove that $\phi_{mD}$, which is a priori only
a rational map, is a morphism). The Iitaka dimension of a divisor is
never greater than the numerical Iitaka dimension~\cite{ev92}, so our
aim is to prove the reverse inequality in the case that $D$ is nef.

Because of Matsushita's Theorem~\ref{matsushita} cited above, a
non-trivial fibration of $X$ has to be over an $n$-dimensional
base. So as soon as some $\H^0(X,\O(mD))$ has dimension greater than
one, and we show $\phi_{mD}$ is a genuine morphism, we can conclude
that the Iitaka dimension is indeed $n$. One way to show a line bundle
has sections is to show that its higher cohomology vanishes, and we
will look at vanishing theorems in the next subsection. First let us
indicate why Conjecture~\ref{nef_sqzero} is important for the study of
deformation classes of holomorphic symplectic manifolds.

\begin{prop}
Every irreducible holomorphic symplectic manifold $X^{\prime}$ with
second Betti number $b_2(X^{\prime})\geq 5$ can be deformed to
another, $X$, such that $X$ contains a non-trivial nef divisor $D$
with $q_X(D)=0$.
\end{prop}

\begin{proof}
The idea is essentially the same as the one used to prove
Proposition~\ref{deform_K3}. The lattice
$(\H^2(X^{\prime},\Z),q_{X^{\prime}})\cong (\Gamma,q_{\Gamma})$ has
rank at least five and is indefinite (more precisely, it has signature 
$(3,b_2(X^{\prime})-3)$). By Meyer's Theorem (see~\cite{cassels78},
page 75) it is isotropic, ie.\ it contains a non-trivial
element whose square is zero. Using an argument similar to those in
Section 4.7 of~\cite{huybrechts02}, we deform to a non-projective $X$
whose Picard group is generated by a primitive divisor $D$ with
$q_X(D)=0$. Either $D$ or $-D$ lies on the boundary of the closure of
the positive cone, which coincides with the closure of the K{\"a}hler
cone by Huybrechts' result~\cite{huybrechts99ii}. Thus either $D$ or
$-D$ is nef, and we are done.
\end{proof}

\begin{rem}
The author is grateful to the referee for pointing out Meyer's
Theorem. There do exist quadratic forms of signature $(3,0)$ and
$(3,1)$ which are not isotropic, and therefore the proof fails when
$b_2(X^{\prime})$ is three or four. We don't currently know whether
there exist irreducible holomorphic symplectic manifolds with those
Betti numbers.
\end{rem}

\subsection{Vanishing theorems}

Huybrechts~\cite{huybrechts99i} proved a Riemann-Roch formula for
holomorphic symplectic manifolds which says
$$\chi(\O(D))=\sum_{i=0}^n a_{2i}(q_X(D))^i$$
for constants $a_{2i}$. Thus given $D$ with $q_X(D)=0$ we have
$$\chi(\O(mD))=a_0=\chi(\O_X)=n+1.$$
Let $\h^i(mD):={\mathrm{dim}}\H^i(X,\O(mD))$; then
$$\sum_{i=0}^{2n} (-1)^i\h^i(mD)=n+1.$$
In the case that $m=1$ and $D$ is non-trivial and primitive, we want
to show that $\h^0(D)=n+1$, and obviously it suffices to show that
$\h^i(D)=0$ for $i>0$ (the calculations in the next subsection
indicate that this is the expected behaviour). In fact, it is enough
to show this for even $i$, as this would imply $\h^0(D)\geq n+1$. So
we have to prove a vanishing theorem for a line bundle with some
special properties (nef and $q_X(D)=0$) on a holomorphic symplectic
manifold.

If $D$ were ample the vanishing result would be immediate. Instead we
are in the `semi-positive' situation of having a nef line bundle. For
projective $X$, the Kawamata-Viehweg vanishing theorem would imply
that $\h^i(D)=0$ for $i$ greater than the numerical Iitaka dimension
of $D$, which in this case is $n$ (for example, see~\cite{ev92}). For
a K3 surface this suffices, as it proves $\h^i(D)$ vanishes for all
even $i>0$; but as soon as $n\geq 2$ we still have some work to
do. This is as much as we can expect to learn from general vanishing
theorems, as we have used both the fact that $D$ is nef and that it
has numerical Iitaka dimension $n$. In the non-projective case, we
even lack a full analogue of the Kawamata-Viehweg theorem at the
moment (cf.\ \cite{dp02}).

To make further progress we need to use some other property of $X$,
and there is really only one possibility, the holomorphic symplectic
form
$$\s\in\H^0(X,\La^2T^*)=\H^{2,0}(X).$$
Take its complex conjugate
$$\bar{\s}\in\H^{0,2}(X)=\H^2(X,\O_X).$$
Verbitsky~\cite{verbitsky90} showed that tensoring with $\s$ gives a
map on cohomology, which together with the adjoint map 
generates a holomorphic Lefschetz action on the cohomology of a
holomorphic symplectic manifold. (This is just part of a more
elaborate ${\mathrm{SO}}(4,1)$ action; see~\cite{verbitsky90}.) Might
it be possible to extend this action in some way to the cohomology of
line bundles on $X$? 

It is certainly possible to define a map on smooth forms with values
in a sheaf $\mathcal{F}$
$$L:{\mathrm C}^{\infty}(\E^{p,q}\otimes\mathcal{F})\rightarrow{\mathrm
C}^{\infty}(\E^{p,q+2}\otimes\mathcal{F})$$
by multiplying with $\bar{\s}$; it also has an adjoint 
$$L^*:{\mathrm C}^{\infty}(\E^{p,q+2}\otimes\mathcal{F})\rightarrow{\mathrm
C}^{\infty}(\E^{p,q}\otimes\mathcal{F}).$$ 
Moreover, $\bar{\s}$ is $\bar{\partial}$-closed (in fact, parallel) so
we get induced maps on cohomology for which we will use the same
notation.

\begin{defn}
Define the {\em sheaf Lefschetz action\/} to be the map on sheaf
cohomology induced by multiplication by $\bar{\s}$
$$L:\H^q(X,\Omega^p\otimes\mathcal{F})\rightarrow\H^{q+2}(X,\Omega^p\otimes\mathcal{F})$$
and the map induced by the adjoint to multiplication by $\bar{\s}$
$$L^*:\H^{q+2}(X,\Omega^p\otimes\mathcal{F})\rightarrow\H^q(X,\Omega^p\otimes\mathcal{F}).$$
\end{defn}

\begin{rem}
A priori, these maps may simply be trivial at the level of
cohomology. 
\end{rem}

We are interested in the case $p=0$ and $\mathcal{F}=\O(D)$ for a nef
divisor $D$. The most desirable behaviour for our purposes is that
$L^*$ is injective. If this is the case, then we see immediately that
$\h^0(D)$ must be at least as large as any other $\h^i(D)$ for $i$
even. For $X$ projective, the Kawamata-Viehweg vanishing theorem gives
$\h^i(D)=0$ for $i>n$. We really only need $\h^{2n}(D)=0$, and in fact
this holds for $X$ non-projective too, provided it is compact and
K{\"a}hler (see Demailly and Peternell~\cite{dp02}). We already know 
$$\sum_{i\mbox{ even}}h^i(D)=n+1+\sum_{i\mbox{ odd}}h^i(D)\geq n+1$$
and since $\h^{2n}(D)=0$ there are at most $n$ non-zero terms on the
left hand side. Therefore $\h^0(D)$ is at least two, and the map
$\phi_D$ has positive dimensional image. As mentioned earlier,
provided $\phi_D$ is a genuine morphism, Matsushita
Theorem~\ref{matsushita} then implies that $D$ has Iitaka dimension
$n$.

\subsection{Some calculations}

Let us compare our `hoped-for' behaviour of the previous subsection
with what happens in practice. Let $\pi:X\rightarrow\P^n$ be an
abelian fibred irreducible holomorphic symplectic manifold, and let
$D$ be the nef divisor with $q_X(D)=0$ given by pulling-back a
hyperplane. We will compute $\h^i(mD)$ using the Leray spectral
sequence, which gives
$$\H^p(\P^n,R^q\pi_*\O(mD))\Rightarrow\H^{p+q}(X,\O(mD)).$$
The left hand side equals
$$\H^p(\P^n,R^q\pi_*\pi^*\O_{\P^n}(m))=\H^p(\P^n,\O_{\P^n}(m)\otimes
R^q\pi_*\O_X).$$
For $q=0$
$$R^0\pi_*\O_X\cong\O_{\P^n}.$$
For $q=1$, observe that the fibre of $R^1\pi_*\O_X$ over a point
$x\in\P^n$ is
$$\H^1(\pi^{-1}(x),\O_{\pi^{-1}(x)})\cong\overline{\H^0(\pi^{-1}(x),\Omega^1_{\pi^{-1}(x)})}.$$
The holomorphic symplectic form on $X$ identifies the tangent $T_X$
and cotangent $\Omega^1_X$ bundles. Restricting to $\pi^{-1}(x)$ these
decompose into 
$$\pi^*T_{\P^n}\oplus
T_{\pi^{-1}(x)}\qquad\mbox{and}\qquad\pi^*\Omega^1_{\P^n}\oplus\Omega^1_{\pi^{-1}(x)}$$
and since the fibres are holomorphic Lagrangian, we can identify
$\Omega^1_{\pi^{-1}(x)}$ with $\pi^*T_{\P^n}$. Thus
$$(R^1\pi_*\O_X)_x\cong\overline{\H^0(\pi^{-1}(x),\pi^*T_{\P^n})}=\overline{(T_{\P^n})}_x.$$
Finally, the K{\"a}hler form identifies the antiholomorphic tangent
space with $(\Omega^1_{\P^n})_x$, and so
$$R^1\pi_*\O_X\cong\Omega^1_{\P_n}.$$
More generally, Matsushita~\cite{matsushita00} proved that
$$R^q\pi_*\O_X\cong\Omega^q_{\P_n}.$$
However, since it is more cumbersome to calculate the cohomology of
these sheaves on $\P^n$ for $q\geq 2$, we will restrict to small
values of $n$.

\begin{exm}
When $X$ is an elliptic K3 surface ($n=1$) all terms in the spectral
sequence vanish apart from $(p,q)=(0,0)$ and $(0,1)$, and we find 
$$\h^i(mD)=\left\{\begin{array}{cl}
                 m+1 & i=0 \\
                 m-1 & i=1 \\
                 0   & i=2. \\
                 \end{array}\right.$$
\end{exm}

\begin{exm}
When $n=2$ we can calculate the $q=2$ term
\begin{eqnarray*}
R^2\pi_*{\mathcal O}_X & \cong &
 R^2\pi_*(\K_X\otimes\pi^*\K_{\P^2}^{\vee}\otimes\pi^*{\mathcal
 O}_{\P^2}(-3)) \\
 & \cong & {\mathcal O}_{\P^2}(-3)\otimes R^2\pi_*\K_{X/ \P^2} \\
 & \cong & {\mathcal O}_{\P^2}(-3)\otimes (R^0\pi_*{\mathcal O}_X)^{\vee}
 \\ 
 & \cong & {\mathcal O}_{\P^2}(-3) 
\end{eqnarray*}
where we have used the triviality of the canonical bundle $\K_X$ and
Serre duality. All terms in the spectral sequence vanish apart from
$(p,q)=(0,0)$, $(0,1)$, and $(0,2)$. We can calculate the dimensions
of these remaining spaces using the Bott-Borel-Weil theorem (for
example, see~\cite{es02}) and hence we find
$$\h^i(mD)=\left\{\begin{array}{cl}
                 (m+2)(m+1)/2 & i=0 \\
                 (m+1)(m-1)   & i=1 \\
                 (m-1)(m-2)/2 & i=2 \\
                 0            & \mbox{otherwise.} \\
                 \end{array}\right.$$
\end{exm}

\begin{exm}
When $n=3$ the calculation is similar to the previous example. We find
$$\h^i(mD)=\left\{\begin{array}{cl}
                 (m+3)(m+2)(m+1)/6 & i=0 \\
                 (m+2)(m+1)(m-1)/2 & i=1 \\
                 (m+1)(m-1)(m-2)/2 & i=2 \\
                 (m-1)(m-2)(m-3)/6 & i=3 \\
                 0            & \mbox{otherwise.} \\
                 \end{array}\right.$$
\end{exm}

\begin{rem}
It is straight-forward to check that $\chi(\O(mD))=n+1$ in each of
these examples. Moreover, for $m=1$ we find $\h^0(D)=n+1$ and
$\h^i(D)=0$ for $i>0$, as expected.
\end{rem}

\begin{rem}
For $n=3$ and large $m$, we see that $\h^2(mD)$ is greater than
$\h^0(mD)$. Therefore the map 
$$L^*:\H^2(X,\O(mD))\rightarrow\H^0(X,\O(mD))$$
cannot be injective. However, for $n=2$ we find $\h^2(mD)<\h^0(mD)$ for
all $m>0$. Moreover, a more detailed analysis seems to indicate that
$L^*$ {\em is\/} injective in this case, suggesting the following
conjecture.
\end{rem}

\begin{conj}
If $E$ is a nef line bundle on an irreducible holomorphic symplectic
$4$-fold $X$, then the map
$$L^*:\H^2(X,E)\rightarrow\H^0(X,E)$$
in the sheaf Lefschetz action is injective.
\end{conj}

As described in the previous subsection, this would go a long way
towards proving Conjecture~\ref{nef_sqzero} for $4$-folds. In the
six-dimensional case, it ought to be possible to analyze the failure
of injectivity of $L^*$, and formulate a slightly weaker statement
instead.

\section{Classification of abelian fibrations which admit sections}

In this section we assume that the irreducible holomorphic symplectic
manifold $X^{2n}$ is fibred by abelian varieties, and that this
fibration has a section. Ultimately, our goal is to relate $X$ to the 
Examples~\ref{HilbK3} and~\ref{genKum} due to
Beauville~\cite{beauville83}.

\subsection{Integrable systems}

Since by Matsushita's Theorem~\ref{matsushita} the fibres of $X$ are
Lagrangian with respect to the holomorphic symplectic form, $X$ is
really an {\em algebraic completely integrable Hamiltonian
system\/}. The local geometry of such systems was studied by
Hurtubise~\cite{hurtubise96}. Under certain hypotheses, he associated
to an integrable system, whose base is an open subset of $\C^n$, a
complex surface $Q$ which admits a holomorphic symplectic form (a
priori, $Q$ is incomplete, though in examples it often extends to an
algebraic surface and the symplectic form to a global meromorphic
form). Moreover, the integrable system is isomorphic (indeed,
symplectomorphic) to the Hilbert scheme of points on $Q$, at least
over a dense open set.

Ultimately, we would like a global version of Hurtubise's result in
the compact case. Namely, we want to associate to our abelian
fibration $X$ a compact complex surface which admits a holomorphic
symplectic form. Such a surface is necessarily a K3 surface or a
torus, and therefore the Hilbert scheme of points on it essentially
gives one of Beauville's Examples~\ref{HilbK3} or~\ref{genKum}. Thus
we wish to relate $X$ to one of Beauville's examples; ideally they
will be birational, and therefore deformation equivalent. However, we
know this cannot be the case with O'Grady's examples, so the
relationship must be more subtle.

Let us briefly review Hurtubise's results~\cite{hurtubise96}. His
first hypotheses is that the abelian fibres are Jacobians of
curves. Later in joint work with Markman~\cite{hm98}, he also dealt
with integrable systems whose fibres are generalized Prym varieties
(note that by a Theorem of Welters, stated as Corollary 2.4 in Chapter
12 of~\cite{bl92}, every principally polarized abelian varieties
occurs as a generalized Prym variety). Let $\J^{2g}\rightarrow
U\subset\C^g$ be a local family of Jacobians of curves, where $U$ is
an open subset of $\C^g$, and let $\S^{g+1}\rightarrow U$ be the
corresponding family of curves. Since $U$ is open in $\C^g$, we can
choose a section of $\S\rightarrow U$, so that each curve has a
base-point. Then the Abel map gives an inclusion
$I:\S\hookrightarrow\J$.

The next hypotheses is that
$$I^*\Omega\wedge I^*\Omega=0$$
on $\S$, where $\Omega$ is the holomorphic symplectic form on
$\J$ (note that $I^*\Omega$ is non-vanishing on $\S$ for dimensional
reasons). This implies that $\S$ is a coisotropic submanifold of $\J$,
and we can quotient by the null foliation on $\S$ to obtain a surface
$Q$. The form $I^*\Omega$ then projects to a holomorphic symplectic
form $\omega$ on $Q$. Hurtubise then proves

\begin{thm}
\label{hurtubise}
There is a birational map
$$\Phi :Q^{[g]}\rightarrow\J$$
which preserves symplectic structures, ie.\ $\Phi$ takes the
holomorphic symplectic form on $Q^{[g]}$ induced by $\omega$ to
$\Omega$ on $\J$. Moreover, for a given curve $S_t$ where $t\in U$,
$\Phi$ takes the Hilbert scheme $S_t^{[g]}={\mathrm{Sym}}^gS_t$ to the
Jacobian $J_t$ of $S_t$ by the Abel map (these are both Lagrangian
submanifolds).
\end{thm}

\begin{rem}
The existence of the birational map $\Phi$ follows from Sklyanin's
separation of variables, and can be viewed as finding Darboux
coordinates on $\J$.
\end{rem}

Hurtubise also studies the geometry of the surface $Q$. He proves that
the curves $S_t$ are embedded in $Q$ and all belong to the same linear
system $|S_t|$. Moreover, if they aren't hyperelliptic curves, then
this linear system embeds $Q$ into $\P^g$ (a similar statement is true
in the hyperelliptic case). For our purposes, the important point is
the following: if a family of Jacobians of curves $\J$ is holomorphic
symplectic, then the base $U$ is an open subset of a linear system of
curves $S_t$ which all lie in a {\em single\/} complex surface $Q$,
and moreover $Q^{[g]}$ is birational to $\J$.

\begin{exm}
An example where these local results extend to global statements is
the Hitchin system~\cite{hitchin87}, the moduli space $\M$ of rank $r$
Higgs bundles over a Riemann surface $\Sig$. Firstly, $\M$ is
isomorphic to the cotangent bundle $T^*\N$ of the moduli space $\N$ of
stable bundles over $\Sig$. Since the latter is a complex manifold,
$T^*\N$ has a canonical holomorphic symplectic structure. As an
abelian fibration, the base of $\M$ is the space of {\em spectral
curves\/}, $r$-to-one finite covers of $\Sig$. If
$\rho:\K_{\Sig}\rightarrow\Sig$ denotes the canonical bundle of
$\Sig$, then we can think of the spectral curves as belonging to the
linear system $|\O(r\rho^*\K_{\Sig})|$ of curves sitting in the total
space $\K_{\Sig}$. In this example, the surface $Q$ is $\K_{\Sig}$,
and there is a birational map between $\M$ and $\K_{\Sig}^{[g]}$,
where $g$ is the genus of the spectral curves. These statements are
global, though the manifolds are all non-compact.
\end{exm}

Another global example comes from the Hilbert scheme of points on a K3
surface $S$.

\begin{exm}
\label{HilbYinZ}
Recall from Example~\ref{abel_fibr_K3} that if $S$ contains a smooth
curve $C$ of genus $g$ then $S^{[g]}$ is birational to an abelian
fibred moduli space $Z$ with base the linear system $|C|\cong\P^g$ and
fibres degree $g$ Picard groups ${\mathrm{Pic}}^g$ of the curves in
the linear system. Associated to the fibration $Z^{2g}\rightarrow\P^g$
is the family of curves $Y^{g+1}\rightarrow\P^g$. The total space $Y$
can be defined by
$$Y:=\{(p,D)|p\in D\in|C|\}$$
in which case the map to $\P^g$ is simply given by $(p,D)\mapsto
D$. Now $Y$ is also foliated, with space of leaves $S$ and leaves
$$Y_q:=\{(p,D)\in Y|p=q\}$$
for $q\in S$. In other words, $Y_q$ is the linear system of curves $D$
which pass through $q$, and hence is isomorphic to $\P^{g-1}$. Thus
$Y$ is both fibred over $\P^g$ and foliated by copies of $\P^{g-1}$.

If $Y$ admits a section, then we can embed $Y$ in $Z$ by mapping
$$(p,D)\mapsto i_*\O_D(p+(g-1)p_0)$$
where $p_0\in D$ is the value of the section at $D\in\P^g$,
$\O_D(p+(g-1)p_0)$ is a degree $g$ line bundle on $D$, and $i$ denotes
the inclusion $D\hookrightarrow S$. In fact, it is enough that $Y$
admit a $(g-1)$-valued multisection. One checks that $\sigma^2|_Y=0$, 
where $\sigma$ is the holomorphic symplectic form on $Z$, and hence
$Y$ is coisotropic. The foliation described above is precisely the
null foliation and quotienting by it allows us to reconstruct the
surface $S$, the space of leaves. Moreover, $Z$ is birational to
$S^{[g]}$. This example is both global and compact.
\end{exm}

\begin{rem}
\label{multisections}
In the local case $\J\rightarrow U$ automatically admits sections. In
the global case, assuming the existence of a section gives us a
base-point in each fibre and allows us to identify it {\em
unambiguously\/} with the Jacobian, which is isomorphic to
${\mathrm{Pic}}^0$. However, in Example~\ref{HilbYinZ} the fibres
were actually degree $g$ Picard groups ${\mathrm{Pic}}^g$. Moreover,
we actually needed a section of $Y$, or at least a $(g-1)$-valued
multisection. Of course, a section of $Y$ will induce a section of
$Z$, and a multisection will induce a multisection. Precisely what
kind of section or multisection we need on $Z$ will therefore vary
according to the dimension, and is a delicate question.
\end{rem}

\subsection{Fibrations and foliations}

Returning to our general abelian fibration $X\rightarrow\P^n$, let us
try to imitate the geometry of Example~\ref{HilbYinZ}. Thus we want
to regard the base $\P^n$ as being a linear system of curves on some
compact complex surface. Let us assume that the fibres of
$X\rightarrow\P^n$ are (isomorphic to) Jacobians of curves, and denote
the family of curves by $Y^{n+1}\rightarrow\P^n$. Moreover, let us
assume that $Y$ can be embedded in $X$ via the Abel map on each fibre
(cf.\ Remark~\ref{multisections}), and that $\sigma^2|_Y$
vanishes. Then $Y$ is coisotropic and we obtain a surface $Q$ by 
quotienting by the null foliation. Moreover, each curve $Y_t$, for
$t\in\P^n$, should project isomorphically to a curve in $Q$, and these
curves should all lie in the same linear system. For $q\in Q$, the
leaf $Y_q\cong\P^{n-1}$ will then correspond to the linear subsystem
of curves passing through $q$.

Thus $Y$ is both fibred by curves over $\P^n$, and foliated by copies
of $\P^{n-1}$ with space of leaves $Q$. The compact surface $Q$ will
admit a holomorphic symplectic form, and ideally we'd like to show it
is a K3 surface. 

\begin{conj}
\label{pp_K3}
Suppose the irreducible holomorphic symplectic manifold $X^{2n}$ is an
abelian fibration, with fibres isomorphic to Jacobians of curves. If
the fibration admits the appropriate multisection (cf.\
Remark~\ref{multisections}) then $X$ is birational to the Hilbert
scheme $S^{[n]}$ of $n$ points on a K3 surface $S$.
\end{conj}

\begin{rem}
When $n=2$ the conjecture has already been proved by
Markushevich~\cite{markushevich96}, and we shall discuss his result in
the next subsection. Note that in this dimension, the condition
$\sigma^2|_Y=0$ is automatic. When $n=5$ O'Grady's space
$\widetilde{\M}$ will possibly be a counter-example to the conjecture,
but recall (Example~\ref{abel_fibr_ogrady}) that we didn't establish a
genuine fibration on $\widetilde{\M}$. Moreover, the fibration may not
admit the appropriate multisection. We will look more closely at the
fibration on $\widetilde{\M}$ later, in Example~\ref{multiple}.
\end{rem}

We saw in Remark~\ref{polarizations} that the generalized Kummer
variety $K_n$ can be deformed to an abelian fibration whose fibres
have polarization of type $(1,\ldots,1,n+1)$; presumably principally
polarized fibres cannot occur. In the projective case, the following
observation is apparently due to Mukai; the author is grateful to
Kieran O'Grady for pointing this out.

\begin{prop}
Suppose the projective abelian fibration $X\rightarrow\P^n$ is a
deformation of the generalized Kummer variety $K_n$. The induced
polarization of the fibres of $X$ cannot be principal.
\end{prop}

\begin{proof}
Let $E$ be the hyperplane section of $X\subset\P^N$, and let $D$ be
the pull-back of a hyperplane from the base $\P^n$. Applying Fujiki's
formula~(\ref{fujiki}) to $E+tD$ gives
\begin{eqnarray*}
\int_X(E+tD)^{2n} & = & c_Xq_X(E+tD)^n \\
                  & = & c_X(q_X(E)+2tq_X(E,D))^n \\
\end{eqnarray*}
and taking the coefficient of $t^n$ we find
$$\int_XE^nD^n=\frac{2^n(n!)^2}{(2n)!}c_Xq_X(E,D)^n.$$
Since $D^n$ is a fibre $F$, the left hand side is just
$\int_F(E|_F)^n$. On the right hand side, the Fujiki constant $c_X$ is
a deformation invariant and for $K_n$ is
$$c_{K_n}=\frac{(n+1)(2n)!}{2^nn!}$$ 
(this is computed in~\cite{ogrady01}, for example). Therefore
$$\int_F(E|_F)^n=(n+1)!q_X(E,D)^n$$
where $q_X(E,D)^n\in\Z$, but if $E|_F$ were a principal polarization
of $F$, the left hand side would equal $n!$.
\end{proof}

\begin{rem}
The polarization of the fibres in Example~\ref{abel_fibr_genK} are
therefore `minimal' according to the formula above (ie.\ correspond to
$q_X(E,D)=1$). By comparison, the Fujiki constant for $S^{[n]}$ is
$$c_{S^{[n]}}=\frac{(2n)!}{2^nn!}$$
and hence the induced polarization of the fibres is principal when
$q_X(E,D)=1$.
\end{rem}

The analogue of Conjecture~\ref{pp_K3} for generalized Kummer
varieties would be the following.

\begin{conj}
Suppose the irreducible holomorphic symplectic manifold $X^{2n}$ is an
abelian fibration, and the fibres have polarization
$(1,\ldots,1,n+1)$. If the fibration admits the appropriate
multisection then $X$ is birational to a generalized Kummer variety
$K_n$.
\end{conj}

The goal would be to again show that the base $\P^n$ of the fibration
is a linear system of curves on a surface (more precisely, a complex
tori $T$), but a priori we don't even have a family of curves. In
Example~\ref{abel_fibr_genK} the fibres were abelian subvarieties
of Jacobians, and perhaps this is one way to approach the problem.

If our goal is to classify holomorphic symplectic manifolds we should
investigate other polarizations too. Presumably there will be some
related to $(1,\ldots,1)$ and $(1,\ldots,1,n+1)$ (for example,
multiples) which will occur for $S^{[n]}$ and $K_n$ respectively, but
the other polarizations will not arise. The justification would be
that if one can associate a holomorphic symplectic surface to the
abelian fibration $X$, then the K3 surface and complex tori already
exhaust all possibilities.

\subsection{Four-folds fibred by Jacobians}

The following theorem is due to Markushevich~\cite{markushevich96}.

\begin{thm}
Suppose the irreducible holomorphic symplectic four-fold
$\pi :X\rightarrow\P^2$ is fibred by Jacobians of genus-two curves,
and that the fibration admits a section. Then $X$ is birational to
$S^{[2]}$ for some K3 surface $S$.
\end{thm}

\begin{proof}
A careful proof may be found in~\cite{markushevich96}. Here let us
just outline how the K3 surface $S$ is constructed: it is actually the
double cover of the dual plane $(\P^2)^{\vee}$ branched over a sextic
$B$.

Let $Y^3\rightarrow\P^2$ be the family of genus-two curves. Each curve
$Y_t$, for $t\in\P^2$, is hyperelliptic, being a
double cover of $\P^1_t:=\P(\H^0(Y_t,\K_{Y_t}))$ branched over six
points. We claim that each of these lines $\P^1_t$ is canonically
embedded in the dual plane $(\P^2)^{\vee}$.

The fibre $X_t$ is the Jacobian of $Y_t$, and therefore its tangent
space $T_pX_t$ at any point $p\in X_t$ is $\H^0(Y_t,\K_{Y_t})$. Using
the holomorphic symplectic form, and the fact that the fibres of $X$
are Lagrangian, we can identify $T_pX_t$ with
$$(\pi^*\Omega^1_{\P^2})_p=(\Omega^1_{\P_2})_t.$$
The Euler sequence on $\P^2$ gives
$$0\rightarrow\Omega^1_{\P^2}\rightarrow\H^0(\P^2,\O(1))\otimes
\O(-1)\rightarrow\O_{\P^2}\rightarrow 0$$
and projectivizing we get the inclusion
$$\P(\Omega^1_{\P^2})\hookrightarrow\P(\H^0(\P^2,\O(1)))=(\P^2)^{\vee}$$
where the right hand side is the trivial bundle over $\P^2$ with fibre
$(\P^2)^{\vee}$. Taking the fibre over $t\in\P^2$ proves the claim.

The six branch points on $\P^1_t$ will vary holomorphically with $t$,
so for a pencil of these lines in $(\P^2)^{\vee}$ it will cut out a
sextic. A priori, different pencils could give different sextics;
however, the six branch points on $\P^1_t$ will actually be the points
of intersection of $\P^1_t$ with the curve $B$ dual to the the
degeneracy locus of the fibration $X$ (ie.\ the fibres of $X$ will
degenerate over some curve $\Delta\in\P^2$, and $B\in(\P^2)^{\vee}$
will be dual to $\Delta$). This establishes that $B$ is a sextic, and
the double cover of $(\P^2)^{\vee}$ branched over $B$ is therefore a
K3 surface $S$. Moreover, pulling-back the line $\P^1_t$ from
$(\P^2)^{\vee}$ will give us a curve in $S$ isomorphic to $Y_t$, as
both curves are double covers of $\P^1_t$ branched over the same six
points.

We have thus realized the base $\P^2$ as a linear system of curves on
a K3 surface. Indeed $X$ is isomorphic to the moduli space $Z$ in
Example~\ref{abel_fibr_K3}, which moreover is birational to
$S^{[2]}$. (One should be careful about singular fibres:
see~\cite{markushevich96} for the details.)
\end{proof}

\begin{rem}
The fact that the curves $Y_t$ are hyperelliptic is no longer true for
higher genus, in general. However, genus-three curves can be realized
as plane quartics, and there is hope in the six-fold case of
constructing a K3 surface as a quartic in $(\P^3)^{\vee}$. On the
other hand, a proof of Markushevich's theorem using Hurtubise's
approach would be more likely to generalization to higher dimensions.
\end{rem}

\section{Abelian fibrations which don't admit sections}

In the previous sections we have looked at the first two steps of our
three-part programme for understanding holomorphic symplectic
manifolds. Our discussion involved some conjectures, some ideas on how
to prove them, and some supporting evidence. The ideas in this section
will be more complete, based on the work of C{\u a}ld{\u
a}raru~\cite{caldararu00}; an example of their application will appear
in the article~\cite{sawon??} currently in preparation. Nevertheless,
there are still some curious phenomena which remain
unexplained (specifically, Example~\ref{multiple}).

\subsection{Gerbes and elliptic fibrations}

Let $X$ be an irreducible holomorphic symplectic manifold which is an
abelian fibration; we don't assume that $X$ admits a section. We wish
to relate to $X$ another abelian fibration which {\em does\/} admit a
section. Ideally we'd like a one-parameter family of deformations
connecting the two manifolds, as in the case of elliptic K3
surfaces. In this subsection we will describe an approach to the
corresponding problem for elliptic fibrations, due to C{\u a}ld{\u
a}raru~\cite{caldararu00}; in the next we will investigate how to
generalize to the case of higher dimensional (abelian) fibres.

Recall from Subsection~\ref{K3_nosection} that the set of all elliptic
K3 surfaces with relative Jacobian $J$ was classified by the Brauer
group $\H^2(J,\O_J^*)$. This group also classifies {\em holomorphic
gerbes\/} on $J$ (see Hitchin~\cite{hitchin01}, for example).

\begin{defn}
Let $\{U_i\}$ be a {\v C}ech cover for $J$ and let $\{\L_{ij}\}$ be a
collection of holomorphic line bundles over $U_{ij}:=U_i\cap U_j$
which satisfy
\ben
\item $\L_{ji}\cong\L_{ij}^{-1}$ for all $i$ and $j$,
\item $\L_{ijk}:=\L_{ij}\otimes\L_{jk}\otimes\L_{ki}$ on $U_{ijk}$
has a trivialization $\beta_{ijk}:U_{ijk}\rightarrow\O^*$ for all $i$,
$j$, and $k$,
\item the trivializations $\{\beta_{ijk}\}$ of $\L_{ijk}$ define a {\v
C}ech 2-cocycle, ie.\ $\delta\beta =1$ (this means that
$\L_{ijkl}:=\L_{jkl}\otimes\L_{ikl}^{-1}\otimes\L_{ijl}\otimes\L_{ijk}^{-1}$
on $U_{ijkl}$ is canonically trivial).
\een
Then $\{\L_{ij}\}$ defines a {\em holomorphic gerbe\/} on $J$. Up to
isomorphism (which we won't define here) the gerbe is classified by
$\beta\in\H^2(J,\O_J^*)$.
\end{defn}

Now suppose $S$ is an elliptic K3 surface with relative Jacobian
$J$. The latter is a moduli space of pure semi-stable sheaves on $S$
(with a certain Hilbert polynomial); we can ask whether there exists a
universal sheaf for this moduli problem. For a single fibre $J_t$, the
Poincar{\'e} line bundle over $S_t\times J_t$ is a universal
bundle, which moreover can be extended over a local fibration. Thus if
$\{U_i\}$ is an affine cover of the base $\P^1$, and
$S_i:=p_S^{-1}(U_i)$ and $J_i:=p_J^{-1}(U_i)$ the corresponding local
fibrations, then there is a local universal bundle $\U_i$ over the
fibre-product $S_i\times_{U_i}J_i$.
$$\begin{array}{ccccc}
   & & {\U}_i & & \\
   & & \downarrow & & \\
   & & S_i\times_{U_i}J_i & & \\
   & \pi_S\swarrow & & \searrow\pi_J & \\
  S_i & & & & J_i \\
   & p_S\searrow & & \swarrow p_J & \\
   & & U_i & & \\
  \end{array}$$
However, these local universal bundles may not patch together to give
a global universal bundle. Above the overlap $U_{ij}$ the local
universal bundles $\U_i$ and $\U_j$ will differ by the pull-back of a
line bundle $\L_{ij}$ from $J_{ij}:=p_J^{-1}(U_{ij})$. One can check
that the collection $\{\L_{ij}\}$ defines a gerbe on $J$, whose
equivalence class we denote by $\alpha_S\in\H^2(J,\O_J^*)$. If
$\alpha_S$ vanishes, the local universal bundles can be patched
together to give a global bundle. Hence $\alpha_S$ is the obstruction
to the existence of a global universal bundle.

To relate this to the classification of elliptic fibrations $S$ we
observe that the choice of a local section $s_i:U_i\rightarrow S_i$
gives a family of basepoints on the fibres $S_t$. This uniquely
defines the Poincar{\'e} line bundle over each fibre, and hence
uniquely defines the local universal bundle $\U_i$. In fact, choosing
a local section is equivalent to choosing a local universal bundle,
and thus our earlier analysis of $S$ based on sections can be
rephrased in terms of universal bundles.

\begin{prop}
The following are equivalent:
\begin{enumerate}
\item $S$ is isomorphic to its relative Jacobian $J$,
\item $S$ admits a section,
\item there exists a universal bundle over $S\times_{\P^1}J$.
\end{enumerate}
The obstruction to all of these existence problems is
$\alpha_S\in\H^2(J,\O_J^*)$.
\end{prop}

\begin{rem}
The approach via gerbes easily generalizes to higher dimensional
elliptic fibrations~\cite{caldararu00}: if the elliptic fibration
$S\rightarrow B$ has relative Jacobian $J\rightarrow B$, then $S$ is
classified by an element of 
$$\H^2(J,\O_J^*)/\H^2(B,\O_B^*).$$
Note that the denominator did not appear in the K3 case as there are
no gerbes on $\P^1$, or indeed on any curve; but $\H^2(B,\O_B^*)$ may
be non-zero as soon as the dimension of $B$ is at least two. In this
case gerbes on $B$ can be pulled-back to $J$, giving an inclusion
$$\H^2(B,\O_B^*)\hookrightarrow\H^2(J,\O_J^*).$$
See C{\u a}ld{\u a}raru's thesis~\cite{caldararu00} for the details.
\end{rem}

\subsection{Gerbes and abelian fibrations}

Now let us consider the case where the fibres themselves are
higher dimensional, namely abelian varieties. We start with a
generalization of the relative Jacobian.

\begin{defn}
Let $X\rightarrow\P^n$ be an abelian fibration. We define the {\em
relative Picard\/} of $X$ to be the moduli space $P$ consisting of
push-forwards of degree-zero rank-one torsion-free sheaves on fibres
of $X$.
\end{defn}

\begin{rem}
In general one has to be careful with singular fibres. Let
$i:X_t\hookrightarrow X$ be the inclusion of a smooth fibre, and let
$i_*L$ be the push-forward of a degree-zero line bundle on $X_t$. Then
more precisely, $P$ should be defined as the moduli space of pure
stable sheaves on $X$ with the same Hilbert polynomial as $i_*L$. If
$P$ is not already complete, then we must add strictly semi-stable
sheaves (which could of course result in a singular space).
\end{rem}

\begin{rem}
On smooth fibres $P$ is given by replacing an abelian variety $X_t$
with its degree-zero Picard group $P_t:={\mathrm{Pic}}^0(X_t)$, ie.\
the dual abelian variety. If $X_t$ is principally polarized then it is
isomorphic to $P_t$ (this is automatically true for elliptic curves),
and moreover $X$ and $P$ will be locally isomorphic as fibrations. For
a general polarization this will no longer be true. A solution to this
problem is to take the relative Picard of $P$: taking the double dual
of an abelian variety returns the original variety. 
\end{rem}

As $P$ is a moduli space of sheaves on $X$, we can investigate the
existence of a universal bundle, like before. For a single (smooth)
fibre $X_t$ and its dual $P_t$, there is a universal bundle on
$X_t\times P_t$. This can be extended to a local universal bundle, at
least away from singular fibres. As in the K3 case, we get a
holomorphic gerbe whose equivalence class $\beta_X\in\H^2(P,\O_P^*)$
is the obstruction to patching these local universal bundles together
to give a global universal bundle. Note that there are no gerbes on
the base $\P^n$, so we needn't worry about $\H^2(\P^n,\O_{\P^n}^*)$.

Denote the relative Picard of $P$ by $X_0$ (we use this notation
because we will actually have $\beta_{X_0}=1$, the trivial gerbe in
$\H^2(P,\O_P^*)$). Then as a fibration, $X$ is locally isomorphic to
$X_0$, at least away from singular fibres, but stuck together in a
different way. Notably, the relative Picard $P\rightarrow\P^n$ has a
canonical section given by the flat family $\O_X$, thought of as the
family of trivial line bundles $\O_{X_t}$ on each fibre of
$X$. Likewise $X_0$ must also have a section $s_0$, as it is the
relative Picard of $P$. Like in Subsection~\ref{K3_nosection},
choosing a local section of $X\rightarrow\P^n$ determines a local
isomorphism with $X_0$, which takes the local section to $s_0$. So $X$ 
is really a torsor over $X_0$, and is thus classified by an element of
$\H^1(\P^n,\B)$ where $\B$ is the sheaf of local holomorphic sections 
of $X_0\rightarrow\P^n$, or equivalently, the sheaf of translations in
fibres.

To relate this to the gerbes above, we observe that local sections of
$X$ are equivalent to local universal bundles. This is almost the same
as the elliptic K3 case, except that a basepoint on an abelian variety
is used to {\em translate\/} a theta divisor, whereas on an elliptic
curve the basepoint {\em is\/} a theta divisor. This essentially
proves the following conjecture: the only remaining concern being the 
singular fibres.

\begin{conj}
The following are equivalent:
\begin{enumerate}
\item $X$ is isomorphic to $X_0$,
\item $X$ admits a section,
\item there exists a universal bundle over $X\times_{\P^n}P$.
\end{enumerate}
The obstruction to all of these existence problems is
$\beta_X\in\H^2(P,\O_P^*)$.
\end{conj}

\begin{rem}
In some examples, $P$ is a smooth irreducible holomorphic symplectic
manifold (in general only smoothness should be in doubt). The
exponential long exact sequence
$$0\rightarrow\H^1(P,\O_P^*)\rightarrow\H^2(P,\Z)\rightarrow\H^2(P,\O_P)\rightarrow\H^2(P,\O_P^*)\rightarrow\H^3(P,\Z)\rightarrow
0$$
and the fact that $\H^2(P,\O_P)\cong\C$ then tells us that
$\H^2(P,\O_P^*)$ will be one-dimensional, with a strange
topology. This is like the K3 case, except that now $\H^3(P,\Z)$ could
be non-trivial and thus $\H^2(P,\O_P^*)$ might not be connected.
\end{rem}

\begin{rem}
The reason that singular fibres cause no problems when constructing
the relative Jacobian $J$ of an elliptic K3 surface $S$ is that
\ben
\item they are isolated,
\item $S$ and $J$ are locally isomorphic away from the singular
fibres.
\een
For elliptic surfaces, the monodromy around the singular fibre
determines it uniquely. For higher dimensional elliptic fibrations the
behaviour of singular fibres is already more complicated, and for
abelian fibrations it is even worse.
\end{rem}

\begin{rem}
The relation between $X_0$ and $P$ is an example of the
Strominger-Yau-Zaslow mirror symmetry
conjecture~\cite{syz96}. Moreover, the relation between $X$ and $P$ is
an example of mirror symmetry with non-trivial $B$-field $\beta_X$. 
Verbitsky~\cite{verbitsky95} showed that generic hyperk{\"a}hler
manifolds (those with Picard number zero) are self-mirror, but our
abelian fibrations have positive Picard number, which is why our
mirror manifolds are not isomorphic in general.

In fact, the main goal of C{\u a}ld{\u a}raru's
thesis~\cite{caldararu00} was to construct {\em twisted Fourier-Mukai
transforms\/}, which combine homological mirror symmetry with the SYZ
conjecture. Our Example~\ref{generic} in the next subsection also
leads to a twisted Fourier-Mukai transform, as will be proved
elsewhere~\cite{sawon??}.
\end{rem}

\subsection{Two examples}

This first example will be described in detail in the
paper~\cite{sawon??} currently in preparation. It has also been
studied by Markushevich~\cite{markushevich96} from a slightly
different point of view.

\begin{exm}
\label{generic}
If the K3 surface contains a smooth genus $g$ curve $C$, then we saw 
in Example~\ref{abel_fibr_K3} that the Hilbert scheme $S^{[g]}$ is
birational to a moduli space $Z$ which is abelian fibred over the base
$|C|\cong\P^g$. The fibre ${\mathrm{Pic}}^gD$ of $Z$ over $D\in|C|$ is
principally polarized, and hence $Z$ is locally isomorphic to its
relative Picard $P$, at least away from singular fibres. 

In the case $g=2$, $S$ is the double cover of the dual plane
$(\P^2)^{\vee}$, branched over a sextic; the curves $D$ are pull-backs
of lines in $(\P^2)^{\vee}$. Write $Z_k$ for the fibration over
$|C|\cong\P^2$ with fibres degree $k$ Picard groups
${\mathrm{Pic}}^k$. Then $Z_k$ is just the Mukai moduli space of
stable sheaves on $S$ with Mukai vector $v=(0,[C],k-1)$. In the
generic situation, $v$ is primitive and hence $Z_k$ is smooth and
compact for all values of $k$. `Generic' here means that the sextic
defining $S$ does not admit a tritangent, for if it admitted a
tritangent then some divisor in $|C|$ would split into two homologous
parts; then $v$ would not be primitive for $k$ odd. It follows also
that $Z_k$ is deformation equivalent to $S^{[2]}$.

The fibrations $Z_k$ are all locally isomorphic to each other. We also
know that $Z_0$ and $Z_2$ admit (canonical) sections, given by taking
$\O_D\in{\mathrm{Pic}}^0$ and $\K_D\in{\mathrm{Pic}}^2$, respectively,
for each curve $D\in|C|$. Thus $Z=Z_2$ is isomorphic to $Z_0$.

We claim that $Z_1$ does not admit a section in the generic case. Thus
$Z_1$ is not isomorphic to $Z_0$ and $Z_2$. In fact we prove this in
the opposite order: using O'Grady's description of the weight-two
Hodge structure of a moduli space of sheaves on a K3
surface~\cite{ogrady97} we can show that $Z_1$ is generically not
isomorphic to $Z_0$ and $Z_2$. Therefore it cannot admit a section. 

All of the fibrations $Z_k$ have the same relative Picard $P$, which
has fibre
$${\mathrm{Pic}}^0({\mathrm{Pic}}^kD)\cong{\mathrm{Pic}}^0D$$
over $D\in|C|$ (the isomorphism is canonical). This shows that $P$ is
isomorphic to $Z_0$ away from the singular fibres, but then $P$ can be
completed by adding the same singular fibres as in $Z_0$, and hence is
deformation equivalent to $S^{[2]}$. The exponential long exact
sequence then shows that $\H^2(P,\O_P^*)$ is one-dimensional and
connected. So although $Z_0$ and $Z_1$ are generically not isomorphic,
they do sit in the same one-dimensional family of deformations (of
course, they are both deformations of $S^{[2]}$). Furthermore, it can
be shown that $Z_1$ admits a 2-valued multisection and hence the gerbe 
$\beta_1\in\H^2(P,\O_P^*)$ classifying $Z_1$ is 2-torsion.
\end{exm}

\begin{rem}
The moduli space of K3 surfaces $S$ which are double covers of the
plane is 19-dimensional. Therefore the manifold $Z=Z_0$, which is the
Mukai flop of $S^{[2]}$, belongs to a 19-dimensional family of abelian 
fibrations, which all admit sections. Then $\H^2(P,\O_P^*)$ adds an
extra dimension, giving us a 20-dimensional family of abelian
fibrations, which sits inside the 21-dimensional moduli space
of deformations of $S^{[2]}$. This is in agreement with
Conjectures~\ref{nef_sqzero} and~\ref{just_sqzero}, which imply that
being an abelian fibration is a codimension one property.

This 20-dimensional family of abelian fibrations also intersects
transversely the set of projective varieties in the moduli space of
$S^{[2]}$ (projectivity is also codimension one). This is because the
projective varieties correspond to torsion points in $\H^2(P,\O_P^*)$,
as in the K3 case (cf.\ Remark~\ref{torsion}).
\end{rem}

Example~\ref{generic} exhibits the expected behaviour: namely, {\em
twisting\/} the fibration $Z_0$ by the non-trivial gerbe $\beta_1$
gives a new fibration $Z_1$ which is not isomorphic to, but is 
nonetheless a deformation of $Z_0$. This is because the singular
fibres were well controlled, which is not always the case. The most
curious effects appear to be produced by {\em multiple
fibres\/}. Recall that elliptic K3 surfaces cannot contain multiple
fibres, so the behaviour in the next example is purely a higher
dimensional phenomena.

\begin{exm}
\label{multiple}
Assume again that the K3 surface $S$ is the double cover of $\P^2$
branched over a generic sextic, but this time let $C$ be the pull-back
of a conic from $\P^2$. Then $C$ is a genus five curve. Denote by
$Z_5$ and $Z_6$ the Mukai moduli spaces $\M^s(0,[C],1)$ and
$\M^{ss}(0,[C],2)$ respectively. The latter has non-primitive Mukai
vector, so we have added strictly semi-stable sheaves to obtain a
singular space. Both $Z_5$ and $Z_6$ are abelian fibrations over
$|C|\cong\P^5$, with fibres ${\mathrm{Pic}}^5D$ and
${\mathrm{Pic}}^6D$ over $D\in|C|$, respectively. Therefore over
smooth curves $D$, $Z_5$ and $Z_6$ are locally isomorphic as
fibrations; one can also show that $Z_6$ admits a section, whereas
$Z_5$ generically does not. 

We expect that $Z_5$ will be obtained from $Z_6$ by twisting, as in
Example~\ref{generic}. However, this cannot be the case because the
singular fibres of $Z_5$ and $Z_6$ are clearly different, and indeed
$Z_6$ is not even smooth. In particular, the pull-back from $\P^2$ of 
a double line gives a non-reduced divisor $2E\in|C|$, and the fibres
$F_5\subset Z_5$ and $F_6\subset Z_6$ over $2E$ are somewhat like
multiple fibres in the theory of elliptic surfaces. In this case,
$F_5$ is certainly not isomorphic to $F_6$ (over simpler kinds of
non-smooth divisors in $|C|$, the singular fibres of $Z_5$ and $Z_6$
might still be isomorphic). If we write $U_5\subset Z_5$ and
$U_6\subset Z_6$ for the unions of the smooth fibres, then we can
perhaps regard $U_5$ as a twisting of $U_6$ (after extending the
theory to incomplete manifolds).

The importance of this example comes from the fact that $Z_5$ is
birational to $S^{[5]}$, and an open subset of $Z_6$ is birational to
an open subset of O'Grady's moduli space $\widetilde{M}$, as we saw in
Examples~\ref{abel_fibr_K3} and~\ref{abel_fibr_ogrady}
respectively. Thus $S^{[5]}$ can be deformed to $Z_5$, and a
symplectic desingularization $\hat{Z}_6$ of $Z_6$ (if it exists) can
be deformed to $\widetilde{M}$; but $S^{[5]}$ and $\widetilde{M}$ are
not deformation equivalent (they have different second Betti
numbers). So the combination of twisting $U_6$ to get $U_5$, and
gluing in different singular fibres, has produced a topological
change. In fact, this is more likely to be due to the different
singular fibres, as the twisting should belong to a one-dimensional
family of deformations.
\end{exm}

\begin{rem}
These two examples are related to a point raised by Beauville
in~\cite{beauville99}. He describes the birational correspondence
between $S^{[g]}$ and $Z_g:=Z$ (cf.\ Example~\ref{abel_fibr_K3}), and
remarks that not much seems to be known about the birational geometry
of the other abelian fibrations $Z_k$, whose fibres over $D\in|C|$ are
${\mathrm{Pic}}^kD$.

In Example~\ref{generic} we saw that $Z_k$ is isomorphic to either
$Z_0$ or $Z_1$ depending on whether $k$ is even or odd. Moreover,
$Z_0$ and $Z_1$ are deformation equivalent but generically not
isomorphic. To prove they are not isomorphic, we used O'Grady's
result~\cite{ogrady97} to show that they have different weight-two
Hodge structures. This also implies they are not birational in the
generic situation, as birational holomorphic symplectic manifolds
represent non-separated points~\cite{huybrechts99i} in the moduli
space and would therefore have isomorphic periods.

In Example~\ref{multiple} we actually find that $Z_5$ and $Z_6$ (after
the appropriate desingularization) are not even deformation
equivalent, and therefore certainly not birational.
\end{rem}

\vspace*{5mm}

{\bf Acknowledgements:} I would like to express my gratitude to the
organizers for the invitation to speak at G{\"o}kova, and to
T{\"U}B{\.I}TAK for its continued support of the conference. I'm
grateful to Daniel Huybrechts, Kieran O'Grady, and Miles Reid for many
useful conversations and helpful suggestions. The referee's comments
were also very helpful.

Until September 2002 the author was a fellow of New College, Oxford,
though some of this work was carried out during visits to Rome and
Cambridge, supported by the European Contract Human Potential
Programme, Research Training Network HPRN-CT-2000-00101 and the Isaac
Newton Institute respectively. Thanks also to Massimiliano Pontecorvo,
the Universit{\'a} di Roma ``La Sapienza'', and King's College,
Cambridge, for hospitality.

\end{document}